\documentclass[a4paper,11pt]{article}
\usepackage[T1]{fontenc}

\date{}
\usepackage{amsmath}
\usepackage{graphicx,color}
\usepackage{makeidx}         
\usepackage{graphicx}
\usepackage{epsfig}
\usepackage[english]{babel}
\usepackage{graphics}
\usepackage{pgf}
\usepackage{amsmath}
\usepackage{amssymb}
\usepackage{multirow}
\usepackage{graphicx,color}
\usepackage{amsfonts}
\usepackage{upgreek}
\usepackage[section]{placeins}
\def\bd{\buildrel}
\def\Pb{{\bf P}}
\def\Eb{{\bf E}}
\def\Var{{\rm Var}}
\def\endsymbol{$\sqcup\mkern-12mu\sqcap$}
\def\done{\ \endsymbol\medskip}
\numberwithin{equation}{section} 

\begin{document}
\title{Goodness-of-fit tests based on sup-functionals of weighted empirical processes}
\author{Natalia Stepanova
and Tatjana Pavlenko\\
\\
\textit{School of Mathematics and Statistics, Carleton University, Canada} \\
\textit{Department of Mathematics, KTH Royal Institute of Technology, Sweden}}

\maketitle

\begin{abstract}
A large class of goodness-of-fit test statistics based on sup-functionals of weighted empirical processes is proposed and studied.
The weight functions employed are Erd\H{o}s-Feller-Kolmogorov-Petrovski  upper-class functions of a Brownian bridge.
 Based on the result of M. Cs\"{o}rg\H{o}, S. Cs\"{o}rg\H{o}, Horv\'{a}th, and Mason
obtained for this type of test statistics, we provide the asymptotic null distribution theory for the class of tests in hand,
and present an algorithm for tabulating the limit distribution functions under the null hypothesis.
A new family of nonparametric confidence bands is constructed for the true distribution function and it is found to perform very well.
The results obtained, together with a new result on the convergence in distribution of the higher criticism statistic,
introduced by Donoho and Jin, demonstrate the advantage of our approach
over a common approach that utilizes a family of regularly varying weight functions.
Furthermore, we show that, in various subtle problems of
detecting sparse heterogeneous mixtures, the proposed test statistics achieve the detection boundary found by Ingster
and, when distinguishing between the null and alternative hypotheses,
perform optimally adaptively to unknown sparsity and size of the non-null effects.
\end{abstract}



\section{Introduction}
In the context of testing the hypothesis of  goodness-of-fit, the weighted empirical process viewpoint
is rather common and very helpful, provided a suitable weight function is used.  For specific types of alternatives,
certain  classical goodness-of-fit tests, including the Kolmogorov-Smirnov and Anderson-Darling tests, may benefit significantly from using proper
weights. For examples of standard and non-standard weight functions, we refer to
\cite{ADarling,Borovkov,Jager,MR,Nikitin}.

In this paper, we study a class of goodness-of-fit test statistics based on the supremum of weighted
empirical processes. The weight functions $q$ employed are the \textit{Erd\H{o}s-Feller-Kolmogorov-Petrovski
(EFKP) upper-class functions}  of a Brownian bridge. As a new class of test statistics, the empirical processes in EFKP weighted sup-norm metrics appeared
for the first time in the work of M. Cs\"{o}rg\H{o}, S. Cs\"{o}rg\H{o}, Horv\'{a}th, and Mason \cite{CsCsHM}.
In our study, we extend this class by allowing any subinterval $I\subseteq (0,1)$ over which the supremum is taken, thus getting a class
of test statistics indexed by two `parameters', the weight function $q$ and the interval $I$. Having a class of test statistics available
gives more flexibility in selecting particular members of the family to
meet specific needs of practical applications.

Asymptotic theory of the tests based on the empirical distribution function (EDF) is commonly
handled by using empirical process techniques and weak convergence theory on metric spaces. The question of weak convergence of
the weighted EDF-based tests turns out to be rather delicate, and adapting even known convergence results
to newly proposed test statistics is not necessarily straightforward. At the same time,  the convergence
in distribution of weighted empirical processes as in Theorem 4.2.3 of
\cite{CsCsHM} (see also Theorem 26.3 (a) in \cite{DasG})
has not been fully explored and used by the statisticians. For an overview and further results along these lines, we refer to Sections 4.5 and 5.5 of \cite{CsH},
and references therein. The latter theorem of M. Cs\"{o}rg\H{o}, S. Cs\"{o}rg\H{o}, Horv\'{a}th, and Mason motivates and provides the background for the current study.

In Section 2, we consider the EDF-based tests standardized by the EFKP upper-class functions of a Brownian bridge, along with the characterization
of this class from  Section 3 of \cite{CsCsHM}. The choice of EFKP upper-class functions ensures that
the corresponding test statistics take on finite values with probability one. A new result on the convergence in distribution of
the EDF-based tests with a commonly used standard deviation-proportional weight function (see Proposition 2.1) provides
further motivation for employing the EFKP upper-class weight functions.
The EDF-based tests standardized by the EFKP upper-class functions are shown to be consistent against a fixed alternative.

Typically, the utility of the test depends on whether or not one can work out the distribution theory of the corresponding test statistic.
In Section 3,  by adapting the results of Section 4 in \cite{CsCsHM}, we establish the asymptotic null
distribution theory for the suggested tests.
The whole class of test statistics is easily seen to be distribution-free. The latter fact allows us to construct a new family
of nonparametric confidence bands for the true distribution function.
In Section 4, we show that, when applied to the problem of detecting sparse heterogeneous mixtures, the entire class of
test statistics is optimally adaptive to an unknown degree of heterogeneity in the Gaussian and non-Gaussian mixture models.
Due to this fact, in various signal detection and classification problems,  the statistics under consideration may be viewed as competitors
to the so-called \textit{higher criticism} statistic, as defined in \cite{DJ}. Section 5 contains some remarks and comments.
An algorithm for tabulating the limit
distributions of the proposed test statistics is given in Section 6, and the proofs of main results are collected in Section 7.
Section 8 is a brief summary of the study.

Some notations used throughout the paper are as follows.
The symbols ${\bd{\mbox\small{\cal D}}\over =}$ and ${\bd  {\mbox\small{\rm a.s.}}\over =}$ are used
for equality in distribution and almost surely, respectively.
The symbols  ${\bd  {\mbox\small{{\cal D}}}\over \rightarrow}$ and ${\bd  {\mbox\small{{P}}}\over \rightarrow}$
denote convergence in distribution and convergence in probability, respectively.
For a set $A$, $\mathbb{I}(A)$ is the indicator of the set $A$.
The notation $a_n\sim b_n$ means that $\lim_{n\to \infty}a_n/b_n=1$, whereas the notation
$a_n\asymp b_n$ means that $0<\liminf_{n\to\infty} (a_n/b_n)\leq\limsup_{n\to \infty} (a_n/b_n)<\infty$.
We use the symbol $\log a$ for the natural (base $e$) logarithm of the number $a$.

\section{Statement of the problem and motivation}
In this section,  we define the main object of our study, the class of test statistics based on weighted empirical process.
Statistical properties of the corresponding tests depend, to a large extent, on the weight function used.
Faced with a large number of conceivable  weight functions, it is of interest to look at the asymptotic behaviour
 of the resulting test statistics. With this in mind, we obtain an interesting convergence result (see Proposition 2.1) that partially
motivates the current study.

\subsection{Class of test statistics}
Let ${X}_{1},{X}_{2},\ldots $ be a sequence of iid random variables with a continuous
cumulative distribution function (CDF) $F$ on $\mathbb{R}$. Denote by
$\mathbb{F}_{n}(t)=n^{-1} \sum_{i=1}^{n} \mathbb{I} ({X}_{i}\leq t)$, $t\in\mathbb{R},$ the
 EDF based on $X_1,\ldots,X_n$. We are interested in testing the hypothesis of goodness-of-fit
\begin{gather}\label{hnull}
 H_{0}: \, F=F_{0}
\end{gather}
against either  a two-tailed alternative $H_{1}: \, F \neq F_{0}$ or an upper-tailed alternative $H^{\prime}_{1}: \, F > F_{0}$.
Before introducing our class of test statistics,  we shall need  some definitions.

\medskip
\noindent
{\bf Definition 1}: A function $w$ defined on $(0,1)$ will be called {\it strictly  positive} if
$$\inf_{\varepsilon \leq u\leq 1-\varepsilon}w(u)>0, \quad  \textrm{for all} \quad \ 0<\varepsilon<1/2.$$

\noindent
\textbf{Definition 2}: Let $q$ be any strictly positive function defined on $(0,1)$ with the property $q(u)=q(1-u)$ for $u\in(0,1/2)$,
which is  nondecreasing in a neighborhood of zero  and nonincreasing in a neighborhood of one. Such a function will be called
an \textit{Erd\H{o}s-Feller-Kolmogorov-Petrovski (EFKP) upper-class function} of a Brownian bridge
$\{ B(u), 0 \leq u \leq 1\}$, if there exists a constant $0\leq b<\infty$ such that
\begin{gather}\label{fv}
\limsup_{u \to 0} | B(u)|/q(u)  \,\,{\bd  {\mbox\small{\rm a.s.}}\over =}\,\, b.
\end{gather}
An EFKP upper-class function $q$ of a Brownian bridge is called a \textit{Chibisov-O'Reilly function}
if $b=0$ in (\ref{fv}).

Note that an EFKP upper-class function does not need to be continuous. For properties of these functions see Section 3 of \cite{CsCsHM}.
 An important example of an EFKP upper-class function with $0<b<\infty$ in (\ref{fv}) is the function
\begin{gather}\label{qq}
q(u)= \sqrt{u(1-u) \log \log (1/(u(1-u)))}.
\end{gather}
Such a choice of $q$ stems from Khinchine's local law of the iterated logarithm, which says that
\begin{gather}\label{lil}
\limsup_{u\to 0}\frac{W(u)}{\sqrt{u\log\log(1/u)}} \,\,{\bd  {\mbox\small{\rm a.s.}}\over =}\,\,\sqrt{2},
\quad \liminf_{u\to 0}\frac{W(u)}{\sqrt{u\log\log(1/u)}} \,\,{\bd  {\mbox\small{\rm a.s.}}\over =}\,\,-\sqrt{2},
\end{gather}
where $\{W(u),0\leq u<\infty\}$ is a standard Wiener process starting at zero. Indeed, relations (\ref{lil}) imply via the
representation of a Brownian bridge
$\{B(u),0\leq u\leq 1\} \,\,{\bd  {\mbox\small{\cal D}}\over =}\,\,\{W(u)-uW(1),0\leq u\leq 1\}$ that, cf. (\ref{fv}),
$$ \limsup_{u\to 0}\frac{|B(u)|}{\sqrt{u(1-u)\log\log(1/u(1-u))}} \,\,{\bd  {\mbox\small{\rm a.s.}}\over =}\,\,\sqrt{2}.$$
 As an example of a Chibisov-O'Reilly function, we mention here the function
\begin{gather}\label{COR}
q(u)=(u(1-u))^{1/2-\nu}, \quad 0<\nu<1/2.
\end{gather}
In practical applications, we recommend to use the weight function as in (\ref{qq}) since, unlike the Chibisov--O'Reilly function  (\ref{COR}),
it does not involve any parameter that has to be (arbitrarily) chosen by the experimenter.

In order to check whether or not a given strictly positive function $q$ as above is an EFKP upper-class function, the following characterization of upper-class functions
may be used (see Theorem 3.3 in \cite{CsCsHM}). Let $q$ be strictly positive function defined on $(0,1)$ with the property $q(u)=q(1-u)$ for $u\in(0,1/2)$,
which is  nondecreasing in a neighborhood of zero  and nonincreasing in a neighborhood of one.
\textit{Then $q$ is an EFKP upper-class function of a Brownian bridge if and only if
$$I(q,c):=\int_0^{1} (x(1-x))^{-1}\exp\left(-c\left(x(1-x)\right)^{-1}q^2(x)\right)\, dx<\infty$$
for some constant $c>0$ or, equivalently, if and only if  
$$E(q,c):=\int_0^{1} \left(x(1-x)\right)^{-3/2} q(x) \exp\left(-c\left(x(1-x)\right)^{-1}q^2(x)\right)\, dx<\infty,$$
for some constant $c>0$ and $\lim_{x\downarrow 0} q(x)/x^{1/2} =\lim_{x\uparrow 1} q(x)/(1-x)^{1/2}=\infty.$}

The integral $I(q,c)$ appeared in the works of \cite{Chib} and  \cite{OR}. The integral $E(q,c)$
appeared in the works of Kolmogorov, Petrovski, Erd\H{o}s, and Feller (see Section 3 of \cite{CsCsHM} for details).

The family of the EFKP upper-class functions is connected to a commonly used family of regularly varying weight functions, which is defined as follows.

\medskip
\noindent
\textbf{Definition 3}: Let $\delta$ be any strictly positive function defined on $(0,1)$ with the property $\delta(u)=\delta(1-u)$ for $u\in (0,1/2)$,
which is nondecreasing in a neighborhood of zero  and nonincreasing in a neighborhood of one. Such a weight function will be called {\it regularly varying with power $\tau \in (0,1/2]$} if for any $b>0$
$$\lim_{t\to0} \frac{\delta(bt)}{\delta(t)} =b^{\tau}.$$

It is clear that the so-called {\it standard deviation proportional (SDP)} weight function
\begin{gather*}
\delta(t) = \sqrt{t(1-t)}
\end{gather*} is regularly varying with power $\tau=1/2$, whereas the Chibisov-O'Reilly function $ \delta(t) = (t(1-t))^{1/2-\nu}$, $\nu \in (0,1/2)$,  is regularly varying with power $\tau=1/2-\nu$.

When dealing with weighted empirical processes, the advantage of using the family of weight functions as in
Definition 2 over that in Definition 3 will be demonstrated in the next sections.

Now, we are ready to define the test statistics of our interests. Back to testing  $H_{0}:\,  F=F_{0}$ versus  $H_{1}:\,  F \neq F_{0}$ or $H_{1}^{'}:\,  F > F_{0}$,
consider the statistics defined by the formulas
\begin{gather*}
T_{n}(q) =\sup_{0<F_0(t)<1} \frac{\sqrt{n} |\mathbb{F}_{n}(t)-F_{0} (t)|}{q(F_{0}(t))},\quad
T_{n}^{+}(q) =\sup_{0<F_0(t)<1} \frac{\sqrt{n} \left( \mathbb{F}_{n}(t)-F_{0} (t) \right )}{q(F_{0}(t))},
\end{gather*}
where $q$ belongs to the family of the EFKP upper-class functions of a Brownian bridge $\{ B(u), 0 \leq u \leq 1\}$.
As $T_{n}(q)$, in this generality, appeared for the first time in the paper of M. Cs\"{o}rg\H{o}, S. Cs\"{o}rg\H{o}, Horv\'{a}th, and Mason \cite{CsCsHM},
the statistics $T_{n}(q)$ and $T_{n}^{+}(q)$ will be called the \textit{two-sided} and \textit{one-sided} {\it Cs\"org\H{o}-Cs\"org\H{o}-Horv\'{a}th-Mason (CsCsHM) statistics},  respectively.
If $H_{0}$ is true, then by the probability integral transformation, for each $n$,
\begin{gather}\label{ts}
T_{n}(q) \,\,{\bd  {\mbox\small{\mathcal D}}\over =}\,\, \sup_{0<u<1} \frac{\sqrt{n} |\mathbb{U}_{n}(u)-u|}{q(u)},\quad
T_{n}^{+}(q) \,\,{\bd  {\mbox\small{\mathcal D}}\over =}\,\,  \sup_{0<u<1} \frac{\sqrt{n} (\mathbb{U}_{n}(u)-u)}{q(u)},
\end{gather}
where $\mathbb{U}_n(u)=n^{-1} \sum_{i=1}^{n} \mathbb{I} ({U}_{i}\leq u)$ with $U_1,\ldots,U_n$ being iid uniform $U(0,1)$ random variables. The corresponding order statistics will be denoted by $U_{(1)}<\ldots <U_{(n)}$.

The test procedures based on $T_{n}(q)$ and $T_{n}^{+}(q)$ are consistent against the alternatives $H_1:F\neq F_0$ and $H_1^{\prime}:F>F_0$, respectively.
Indeed, consider, for instance, testing $H_0:F=F_0$ versus $H_1:F\neq F_0$ and observe that for any fixed alternative $F\neq F_0$ the statistic
$\|(\mathbb{F}_n-F_0)/q(F_0)\|_{\infty}:=\sup_{0<F_0(t)<1}|\mathbb{F}_n(t)-F_0(t)|/q(F_0(t))$
satisfies
$$\|(\mathbb{F}_n-F_0)/q(F_0)\|_{\infty}\geq \|({F}-F_0)/q(F_0)\|_{\infty}+o_P(1),$$
implying $T_n(q) \,\,{\bd  {\mbox\small{P}}\over \rightarrow}\,\, \infty$ whenever $F\neq F_0$.
From this, 
$$\Pb_{H_1}\left(T_n(q)>t_\alpha(q)  \right)\to 1,\quad n\to \infty.$$
The case of the upper-tailed alternative is treated similarly.

It might be also of interest to consider the following generalization of the  CsCsHM statistics. For $0 \leq a <b \leq 1$, denote $I=(a,b)$ and define the statistics
\begin{eqnarray*}
 T_{n}(q,I)&=& \sup_{a<F_0(t)<b} \frac{\sqrt{n} |\mathbb{F}_{n}(t)-F_{0} (t)|}{q(F_{0}(t))},\\
 T_{n}^{+}(q,I)&=& \sup_{a<F_0(t)<b} \frac{\sqrt{n} (\mathbb{F}_{n}(t)-F_{0} (t))}{q(F_{0}(t))},
 \end{eqnarray*}
which, for each $n$, under the null hypothesis,  have the same distributions as the random variables $\sup_{u\in I} {{\sqrt{n} |\mathbb{U}_{n}(u)-u|}}/{q(u)}$ and
$\sup_{u\in I} {{\sqrt{n} (\mathbb{U}_{n}(u)-u)}}/{q(u)}$, respectively.

\subsection{Asymptotic properties of test statistics with the SDP weight function: connection to the higher criticism approach}

\noindent
The EDF-based tests standardized by the SDP weight function
$\delta(t) = \sqrt{t(1-t)}$ have been extensively studied in the literature (see
\cite{ADarling,Borovkov,DJ,Eicker, J, Jager}, etc.)
A popular statistic of this kind is the higher criticism  statistic, which is defined as
\begin{gather}\label{DJ}
\mbox{HC}_n=\sup_{0<u<\alpha_0} \frac{\sqrt{n} (\mathbb{U}_{n}(u)-u)}{\sqrt{u(1-u)}}, \quad 0<\alpha_0<1.
\end{gather}
The statistic ${\rm HC}_n$ was introduced by Donoho and Jin \cite{DJ} for multiple testing
situations where most of the component problems correspond to
the null hypothesis and there may be a small fraction of component problems that correspond to non-null hypotheses.
The situations of this kind, where there are many independent null hypothesis $H_{0i}$, $i=1,\ldots,n$,
and we are interested in rejecting the joint null hypothesis $\cap_{i=1}^n H_{0i}$,  are considered in Section 4.
Such a multiple testing problem
is closely connected to the testing problem for the Bayesian alternative
studied by Ingster \cite{Ing97}; the latter problem is of importance in various applications
for multi-channel detection and communication systems. It also relates to the problem of optimal feature selection in
a high-dimensional classification framework, as studied in \cite{DJ2009}, where
another version of the higher criticism statistic has been discussed.

The test statistic $\mbox{HC}_n$
is derived from the random variable
$$ \max_{0<\alpha\leq\alpha_0}\frac{\sqrt{n}\left(M_n/n-\alpha\right)}{\sqrt{\alpha(1-\alpha)}},$$
where $M_n$ is the number of hypotheses among $H_{0i}$, $i=1,\ldots, n$, that are rejected at level $\alpha$,
which measures the maximum deviation of the observed proportion of rejections from what one would expect it to be purely by chance
as the Type I error level changes from zero to $\alpha_0$ (see \cite{DJ} and Section 34.7 of \cite{DasG} for details).
Thus, the parameter $\alpha_0$ in (\ref{DJ}) defines a range of significance levels
in multiple-comparison testing and therefore is a number like $0.1$ or $0.2$.

Below we will show that for all large enough $n$ the
distributional properties of $\mbox{HC}_n$ remain the same no matter what $\alpha_{0}$ is chosen.
Therefore, in this context, the  interpretation of $\alpha_{0}$ as a significance level is somewhat misleading.

The convergence properties of the statistic $\mbox{HC}_n$ are  largely determined by the behaviour
of ${\sqrt{n} (\mathbb{U}_{n}(u)-u)}/{\sqrt{u(1-u)}}$ in the vicinity of zero and one:  the latter inflates
when $u$ is close to zero and one.
The almost sure rate at which $\mbox{HC}_n$ blows up is considered in Chapter 16 of \cite{Shorak}.

 To overcome this problem, Donoho and Jin \cite{DJ} suggested to truncate the range over which the supremum in (\ref{DJ})
 is taken to $(1/n,\alpha_0)$;  this resulted in the test statistic
 \begin{gather}\label{DJT}
 \mbox{HC}^{+}_n=\sup_{1/n<u<\alpha_{0}} \frac{\sqrt{n} (\mathbb{U}_{n}(u)-u)}{\sqrt{u(1-u)}}, \quad 0<\alpha_0<1.
\end{gather}
A seemingly better modification of the higher criticism statistic ${\rm HC}_n$ has the form (see, e.g., \cite{Jager})

\begin{gather}\label{DJ2}
\mbox{HC}^*_n=\sup_{U_{(1)}<u<U_{([\alpha_0 n])}} \frac{\sqrt{n} (\mathbb{U}_{n}(u)-u)}{\sqrt{u(1-u)}},\quad 0<\alpha_0<1.
\end{gather}

Unfortunately, truncating the range as in (\ref{DJT}) and (\ref{DJ2}) does not eliminate the problem.
Indeed, all three statistics, $\mbox{HC}_n$,  $\mbox{HC}_n^+$, and  $\mbox{HC}_n^*$, when normalized as
in Eicker \cite{Eicker} and Jaeschke \cite{J},  under the null hypothesis. will have an extreme value distribution
as $n\to \infty$.
The precise statement, due to Eicker and Jaeschke, is as follows (see  Theorem 2 on p. 118 in
\cite{Eicker}, and Theorem on p. 109 in \cite{J}): for any $x\in\mathbb{R}$,
\begin{multline}
\lim_{n\to \infty}\Pb\left(a_n \sup_{0<u<1} \frac{\sqrt{n}(\mathbb{U}_n(u)-u)}{\sqrt{u(1-u)}}-b_n\leq x\right) \\=
\lim_{n\to \infty}\Pb\left(a_n \sup_{U_{(1)}<u<U_{(n)}} \frac{\sqrt{n}(\mathbb{U}_n(u)-u)}{\sqrt{u(1-u)}}-b_n \leq x\right) = E^2(x).\label{j1}
\end{multline}
where
\begin{gather*}
a_n=\sqrt{2 \log\log n},\quad b_n=2\log\log n+\frac12\log\log\log n-\frac12 \log(4\pi),
\end{gather*}
and $E(x)=\exp(-\exp(-x))$ is the
extreme value CDF. It is clear (see also Corollary 2 on p. 110 in \cite{J}) that the above convergence remains valid
when the supremum is taken over $(1/n,1-1/n)$.
We note in passing that Jager and Wellner \cite{Jager} also study appropriately normalized version of $(1/2)({\rm HC}_n^*)^2$ with
${\rm HC}_n^*$ as in (\ref{DJ2}) in terms of an extreme value distribution, cf. their Theorem 3.1.
 Furthermore, the following result holds true (compare to the claim in Section 3 of \cite{DJ}).

\medskip

\noindent \textbf{Proposition 2.1.}
\textit{For any $0<\alpha_0<1$ and any $x\in\mathbb{R}$,
\begin{gather}\label{j2}
\lim_{n\to \infty}\Pb\left(a_n \sup_{0<u<\alpha_0} \frac{\sqrt{n}(\mathbb{U}_n(u)-u)}{\sqrt{u(1-u)}}-b_n \leq x\right) = \exp\left(-\frac12\exp(-x)\right),
\end{gather}
where $a_n$ and $b_n$ are as in {\rm (\ref{j1})}.}

\medskip
\noindent \textit{Remark 2.1.}
{A similar result holds true for the supremum of
$\dfrac{\sqrt{n}|\mathbb{U}_n(u)-u|}{\sqrt{u(1-u)}}$. Namely, with  $a_n$ and $b_n$ as above,
 for any $0<\alpha_0<1$ and any $x\in\mathbb{R}$,
\begin{gather*}
\lim_{n\to \infty}\Pb\left(a_n \sup_{0<u<\alpha_0} \frac{\sqrt{n}|\mathbb{U}_n(u)-u|}{\sqrt{u(1-u)}}-b_n \leq x\right) = \exp\left(-\exp(-x)\right).
\end{gather*}
Based on the result of Darling and Erd\H{o}s {\rm \cite{Darling}}, the proof of this statement exploits the same route as that of Proposition 2.1.}

\medskip

The message from Proposition 2.1 is that, regardless of a particular value of $0<\alpha_0<1$,
one always has the same extreme value distribution on the right side of (\ref{j2}).
Furthermore, as readily seen from the proof of Proposition 2.1, on the left side of (\ref{j2}) either one of the intervals $(1/n,\alpha_0)$ and $(U_{(1)},U_{([\alpha_0 n])})$
can be taken in place of  $(0,\alpha_0)$.
Together with (\ref{j1}) this implies that, in the sup-norm scenario, when normalizing the process ${\sqrt{n}(\mathbb{U}_n(u)-u)}$ by ${\sqrt{u(1-u)}}$, one arrives at
the situation where ``all the action takes place on the tails but, unfortunately, near infinity''.
This phenomenon seems to be well known to the experts.
However, the proof (if any) of statement (\ref{j2}) is not easily accessible. Therefore we shall prove it here in Section 7.

\subsection{Motivation} In a number of recent studies that utilize the higher criticism strategy
(see, for example, \cite{Cai, DJ, DJ2009, Klaus}) the existence of the extreme value approximation
 to the properly normalized higher criticism statistics ${\rm HC}_n$, ${\rm HC}^+_n$, and  ${\rm HC}_n^*$ has, in fact, never been used.
This could be possibly explained by the fact that the convergence in distribution in (\ref{j1}) and (\ref{j2}) is in general slow (see \cite{J} for details and references).
At the same time, without using the extreme value approximation, one immediately faces the problem of how to choose critical values.

 Next, as follows from Proposition 2.1, in the limit, the parameter $0<\alpha_0<1$ that enters
 the higher criticism statistics ${\rm HC}_n$, ${\rm HC}^+_n$, and  ${\rm HC}_n^*$ as the maximum significance level loses its role.
Furthermore, when using the higher criticism statistics in practice, one has to eliminate the influence of the end-point zero by arbitrarily
truncating the interval $(0,\alpha_0)$ over which the supremum in (\ref{DJ}) is taken.

Thus, the application of the higher criticism statistic ${\rm HC}^+_n$ in practice in such a way that its good asymptotic properties,
including optimal adaptivity (see, e.g., Theorem 1.2  in \cite{DJ}),
are preserved is not straightforward.
As noticed in Remark on p. 611 of \cite{DasG}, ``It is not clear for what $n$ the asymptotics start to give reasonably accurate description
of the actual finite sample performance and actual finite sample comparison... Simulations would be informative and even necessary. But the range in which $n$ has to be in order that
the procedure work well when the distance between the null hypothesis and alternative so small
would make the necessary simulations time consuming.''

All this has motivated us to search for a better weighted analog of the higher criticism statistic, for which
the ``action is shifted somewhat to the middle, while properly regulated on the tails''.
With this focus in mind, we propose a different approach, based on the result of  M. Cs\"org\H{o}, S. Cs\"org\H{o}, Horv\'{a}th, and Mason
 on the convergence in distribution  of the uniform empirical process in weighed sup-norm (see Theorem 4.2.3 in \cite{CsCsHM}).
 The test procedures proposed in this paper do not require an unrealistically large sample size of $n=10^6$ and work well even for $n=10^2.$

\section{Asymptotic distribution theory under the null hypothesis}
To control the probability of Type I Error, one typically chooses a test statistic in such a way
that its distribution under the null hypothesis is known, at least, in the limit.
In this section, we present the null limit distributions of $T_n(q)$ and $T_n^+(q)$ and their modifications.

\subsection{Limit distributions of the CsCsHM-type test statistics under the null hypothesis}
The limit distributions of the statistics
$T_{n}(q)$ and $T_{n}^{+}(q)$ under the null hypothesis are easily accessible.
To clarify this claim, we shall need two facts.

\medskip
\noindent \textbf{Fact 1 (Theorem 4.2.3 in \cite{CsCsHM}).} \textit{Let $q$ be a strictly positive function on $(0,1)$
such that it is nondecreasing in a neighbourhood of zero and nonincreasing in a neighbourhood of one. The sequence of random variables
$\sup_{0<u<1}\sqrt{n}|\mathbb{U}_n(u)-u|/q(u)$ converges in distribution to a nondegenerate random variable if and only if
$q$ is an EFKP upper-class function. The latter nondegenerate random variable must be the random variable
$\sup_{0<u<1}|B(u)|/q(u)$.}

\medskip
\noindent \textbf{Fact 2 (Lemma 4.2.2 in \cite{CsCsHM}).} \textit{Whenever $q$ is an EFKP upper-class function, then for each $-\infty<x<\infty$
and any Brownian bridge $B$}
\begin{gather*}
\Pb\left(\sup_{1/n\leq u\leq 1-1/n}|B(u)|/q(u)\leq x \right) \to \Pb\left(\sup_{0<u<1}|B(u)|/q(u)\leq x \right),\quad n\to \infty.
\end{gather*}

\medskip

It now follows from Facts 1 and 2 that, if $H_0$ is true, then as $n\to \infty$  
\begin{eqnarray}
T_{n}(q)&=&\sup_{0<F_0(t)<1} \frac{\sqrt{n} |\mathbb{F}_{n}(t)-F_{0} (t)|}{q(F_{0}(t))} \,\,{\bd  {\mbox\small{{\cal D}}}\over \rightarrow}\,\,  \sup_{0<u<1} |B(u)|/q(u),
\label{mc1}\\
T_{n}^+(q)&=&\sup_{0<F_0(t)<1} \frac{\sqrt{n} (\mathbb{F}_{n}(t)-F_{0} (t))}{q(F_{0}(t))} \,\,{\bd  {\mbox\small{{\cal D}}}\over \rightarrow}\,\, \sup_{0<u<1} B(u)/q(u),
\label{mc2}
\end{eqnarray}

Indeed, in view of (\ref{ts}), the first statement (\ref{mc1}) is just Fact 1 cited above.
Next,  by positivity of the random variable $\sup_{0<u<1}B(u)/q(u)$, Fact 2 continues to hold with
$B(u)/q(u)$ in place of $|B(u)|/q(u)$.
Therefore, having this modification of Fact 2, one of the key elements in the proof of Fact 1, we immediately arrive at
the second statement (\ref{mc2}).
Furthermore, the inspection of the proof of Fact 1 yields that both statements can be generalized by admitting subintervals.
Namely, the following result holds true.

\medskip

\noindent\textbf{Proposition 3.1.}
\textit{Let $q$ be an EFKP upper-class function of a Brownian bridge. Then, under $H_0$, for any numbers $0\leq a<b\leq 1$, as $n\to \infty$,}
\begin{gather*}
\sup_{a<F_0(t)<b} \frac{\sqrt{n}|\mathbb{F}_{n}(t)-F_{0} (t)|}{q(F_{0}(t))} \,\,{\bd  {\mbox\small{{\cal D}}}\over \rightarrow}\,\, \sup_{a<u<b} \frac{|B(u)|}{q(u)}, \nonumber\\
\sup_{a<F_0(t)<b} \frac{\sqrt{n} (\mathbb{F}_{n}(t)-F_{0} (t))}{q(F_{0}(t))} \,\,{\bd  {\mbox\small{{\cal D}}}\over \rightarrow}\,\,
\sup_{a<u<b} \frac{B(u)}{q(u)}.
\end{gather*}

\medskip
For the intervals $(a,b),$ where one of the end-points is either zero or one, with the above mentioned modification of Fact 2 in mind,
both statements are proved exactly along the lines of the proof of Theorem 4.2.3 in \cite{CsCsHM}.
For the interval $(a,b)$ with $0<a<b<1$, the proof is even easier, as Fact 2 is not needed anymore.
Therefore the proof of Proposition 3.1 is omitted.

\medskip

\noindent \textit{Remark 3.1.}
{Define the statistics
\begin{gather}\label{That}
\hat{T}_n(q)=\sup\limits_{0<F_0(t)<1} \dfrac{\sqrt{n}|\mathbb{F}_n(t)-F_0(t)|}{q({\mathbb{F}_n(t)})},
\quad \hat{T}_n^+(q)=\sup\limits_{0<F_0(t)<1} \dfrac{\sqrt{n}(\mathbb{F}_n(t)-F_0(t))}{q({\mathbb{F}_n(t)})},
\end{gather}
where we set ${\sqrt{n}|\mathbb{F}_n(t)-F_0(t)|}/{q({\mathbb{F}_n(t)})}=0$ for $\mathbb{F}_n(t)\in\{0,1\}.$
With a bit of technical work based on the Glivenko-Cantelli theorem and Slutsky's lemma, one can show that
for any \textit{continuous} EFKP upper-class function $q$ on $(0,1)$,  including  the function
$q$ as in {\rm(\ref{qq})}, the statement of Proposition 3.1 remains valid with $\hat{T}_n(q)$ and  $\hat{T}_n^+(q)$ in place of ${T}_n(q)$ and $T_n^+(q)$,
respectively. The latter fact will be used
to construct  a nonparametric confidence band in Section 3.2.}
\medskip

It follows from Proposition 3.1 that the distribution of the empirical process in weighted sup-norms under
consideration depends on the interval over which the supremum is taken, as all the action now
takes place in the middle.
The same observation applies to the two-sided statistic $T_n(q)$.
Therefore, now, when dealing with the statistic
\begin{gather}\label{nasha}
T^+_{n}(q,(0,\alpha_0))=\sup_{0<F_0(t)<\alpha_0} \frac{\sqrt{n}(\mathbb{F}_{n}(t)-F_{0} (t))}{q(F_{0}(t))},\quad 0<\alpha_0<1,
\end{gather}
we may indeed think of the interval $(0,\alpha_0)$ as the range over which significance levels vary in a multiple testing problem, as suggested by
Donoho and Jin \cite{DJ}.

For the normalized empirical process ${\sqrt{n}(\mathbb{F}_n(t)-F_0(t))}/{\sqrt{F_0(t)(1-F_0(t))}}$ the situation is different. In this case, all the action takes place in the tails,
and for any $0<\alpha_0\leq 1/2$ (see Corollary 3 in \cite{J})
\begin{gather*}
a_n \sup_{\alpha_0<u<1-\alpha_0} \frac{\sqrt{n}(\mathbb{U}_n(u)-u)}{\sqrt{u(1-u)}}-b_n  \,\,{\bd  {\mbox\small{{\rm P}}}\over \rightarrow}\,\, -\infty,
\end{gather*}
where $a_n$ and $b_n$ are as in (\ref{j1}).

The convergence results (\ref{mc1}) and (\ref{mc2}) suggest the following test procedures of
asymptotic level $\alpha$. Set $$T(q):=\sup_{0<u<1} |B(u)|/q(u),\quad T^{+}(q):=  \sup_{0<u<1} B(u)/q(u).$$
Then, one would reject $H_{0}$  in favor of $H_1$ when
$T_{n}(q)>t_{\alpha}(q)$, where the critical point $t_{\alpha}(q)$ is chosen to have  $P(T(q) \geq t_{\alpha}(q))=\alpha$; and
one would reject $H_{0}$ in favour of $H^{\prime}_1$  whenever $T^{+}_{n}(q)>t^{+}_{\alpha}(q)$, where
$t^{+}_{\alpha}(q)$ is determined by ${P}(T^{+}(q) \geq t^{+}_{\alpha}(q))=\alpha$.

The main advantage of using the family of test statistics as in (\ref{ts}) is the identification of the limit distribution
under the null hypothesis.
This distribution is tabulated in Appendix.
Although no analytical results on the convergence rates in (\ref{j1}) and (\ref{mc2}) seem to exist,
the confidence bands obtained in the Section 3.2 confirm better convergence properties of  (\ref{mc2})
as compared to  (\ref{j1}). Heuristically, the situation is as follows.
The renormalization of the sup-functional of the standardized uniform empirical process $\{\sqrt{n}(\mathbb{U}_n(u)-u)/\sqrt{u(1-u)}, 0<u<1\}$
in Proposition 2.1 is to pull it back from disappearing to $-\infty$  via obtaining an extreme value
distribution by bounding it away from zero, cf. (\ref{star}).
The thus obtained extreme value distribution, however, appears to be less concentrated ``in the middle''
as compared to that of the CsCsHM statistic; the confidence bounds for the former
tend to be wider ``in the middle'' than those for the latter. On the tails, they appear
to be doing a similar job, with the Eicker-Jaeschke bounds seemingly better there (see Figure~\ref{figure1} on page~\pageref{figure1}).
We note in passing that the Eicker and Jaeschke solution amounts to renormalization
as compared to $n^{1/2}$, while the Cs\"{o}rg\H{o}, Cs\"{o}rg\H{o}, Horv\'{a}th, and Mason approach
is reweighing the standardized uniform empirical process as compared to $(u(1-u))^{-1/2}$.

\subsection{Confidence bands} In this subsection, we compare numerically a new family of confidence bands derived from (\ref{mc1}) with
those of Kolmogorov and Smirnov, and Eicker and Jaeschke.
As before, suppose that $X_1,X_2,\ldots$ is a sequence of i.i.d.
random variables with a common continuous CDF $F(t),$ and let $\mathbb{F}_n(t)=n^{-1}\sum_{i=1}^n \mathbb{I}(X_i\leq t)$ be the EDF.
For the purpose of constructing a confidence band for $F$, we choose the weight function $q$ on $(0,1)$ to be
$$q(u)=\sqrt{u(1-u)\log\log(1/u(1-u))}.$$
 By Remark 3.1, as $n\to \infty$,
\begin{gather*}
\sup\limits_{0<F(t)<1} \dfrac{\sqrt{n}|\mathbb{F}_n(t)-F(t)|}{q({\mathbb{F}_n(t)})}\,\,{\bd  {\mbox\small{{\cal D}}}\over \rightarrow}\,\,  \sup_{0<u<1} |B(u)|/q(u).
\end{gather*}
Therefore, setting
$H(t)=\Pb\left(\sup\limits_{0<u<1}{|B(u)|}/{q(u)}\leq t\right)$
and denoting by $c_{\alpha}$ the $(1-\alpha)$th quantile of $H$, we obtain
\begin{multline*}
1-\alpha=\lim_{n\to \infty}\Pb_F\left(\sup\limits_{0<F(t)<1} \dfrac{\sqrt{n}|\mathbb{F}_n(t)-F(t)|}{q({\mathbb{F}_n(t)})} \leq c_{\alpha}\right)\\
= \lim_{n\to \infty}\Pb_F\left( \mathbb{F}_n(t)-\frac{c_{\alpha}}{\sqrt{n}}\,q( \mathbb{F}_n(t)) \right.\leq F(t) \leq\\
\left.\leq  \mathbb{F}_n(t)+\frac{c_{\alpha}}{\sqrt{n}}\,q( \mathbb{F}_n(t)),\;\forall\; t\in[X_{(1)},X_{(n)} )\right).
\end{multline*}
Since $F(t)$ takes on its values in $[0,1]$, it now follows that
\begin{multline*}
\lim_{n\to \infty}\Pb_F\left(\max\left\{0,\mathbb{F}_n(t)-\frac{c_{\alpha}}{\sqrt{n}}\,q( \mathbb{F}_n(t))\right\}\leq F(t)\leq\right.\\
\left.\leq
\min\left\{1,\mathbb{F}_n(t)+\frac{c_{\alpha}}{\sqrt{n}}\,q( \mathbb{F}_n(t))\right\},\;\forall\;t\in[X_{(1)},X_{(n)} )\right)=1-\alpha.
\end{multline*}
This gives us an asymptotically correct $100(1-\alpha)\%$ confidence band $[L_n(t),U_n(t)]$ for $F(t)$ on the interval $t\in[X_{(1)}, X_{(n)})$, where
$$L_n(t)=\max\left\{0,\mathbb{F}_n(t)-\frac{c_{\alpha}}{\sqrt{n}}q( \mathbb{F}_n(t))\right\},\quad
U_n(t)=\min\left\{1,\mathbb{F}_n(t)+\frac{c_{\alpha}}{\sqrt{n}}q( \mathbb{F}_n(t))\right\},$$
and, given $\alpha\in(0,1)$, the value of $c_\alpha$ is found from Table III in \cite{OP}.
For instance, $c_{0.05}=4.57.$

Now we compare numerically the confidence band obtained with the  $100(1-\alpha)\%$ confidence band $[L_{n,KS}(t), U_{n,KS}(t)]$ for $F(t)$
based on the two-sided Kolmogorov-Smirnov statistic.
We know that
\begin{gather*}
\Pb_F\left(\sqrt{n}\sup_{-\infty<t<\infty}|\mathbb{F}_n(t)-F(t)|\leq x\right)\to K(x), \quad x\in\mathbb{R},
\end{gather*}
where $K(x)=\sum_{k=-\infty}^{\infty}(-1)^k e^{-2k^2 x^2}$ for $x>0$, and zero otherwise, is the Kolmogorov function
which is tabulated in many textbooks on mathematical statistics. Then, the corresponding
lower and upper bounds are given by
$$L_{n,KS}(t)=\max\left\{0,\mathbb{F}_n(t)-\frac{k_{\alpha}}{\sqrt{n}}\right\},\quad
U_{n,KS}(t)=\min\left\{1,\mathbb{F}_n(t)+\frac{k_{\alpha}}{\sqrt{n}}\right\}.$$
Here $k_{\alpha}$ is the $(1-\alpha)$th quantile of the Kolmogorov function $K(x)$. For instance, $k_{0.05}=1.35.$

For all $n$ large enough, we can confine the region where the Kolmogorov-Smirnov confidence band is constructed to the interval $[X_{(1)}, X_{(n)})$.
Indeed, setting
\begin{eqnarray*}
D^{(1)}_{n}&=&\sqrt{n}\sup_{-\infty<t<X_{(1)}}|\mathbb{F}_n(t)-F(t)|,\\
D_n^{(2)}&=&\sqrt{n}\sup_{X_{(1)}\leq t<X_{(n)}}|\mathbb{F}_n(t)-F(t)|,\\
D_n^{(3)}&=&\sqrt{n}\sup_{X_{(n)}\leq t<\infty}|\mathbb{F}_n(t)-F(t)|,
\end{eqnarray*}
we can write $\sqrt{n}\sup_{-\infty<t<\infty}|\mathbb{F}_n(t)-F(t)|=\max\left(D_n^{(1)}, D_n^{(2)}, D_n^{(3)} \right),$
where $D_n^{(1)}\,\,{\bd  {\mbox\small{{\cal D}}}\over =}\,\, \sqrt{n}U_{(1)}$ and $D_n^{(3)}\,\,{\bd  {\mbox\small{{\cal D}}}\over =}\,\,\sqrt{n}(1-U_{(n)})$.
Since $nU_{(1)}$ and $n(1-U_{(n)})$  have exponential limit distributions, it follows that
$D_n^{(1)}$ and $D_n^{(3)}$ converge in probability to zero, and hence
the limit distribution of $\sqrt{n}\sup_{-\infty<t<\infty}|\mathbb{F}_n(t)-F(t)|$ coincides with that of $D_n^{(2)}$.

Next, in order to obtain the Eicker-Jaeschke confidence band, consider the random variable
$$\hat{T}_n=\sup\limits_{0<F(t)<1} \dfrac{\sqrt{n}|\mathbb{F}_n(t)-F(t)|}{\sqrt{\mathbb{F}_n(t)(1-\mathbb{F}_n(t))}},$$
with the said convention that  ${{\sqrt{n}|\mathbb{F}_n(t)-F(t)|}}/{{\sqrt{\mathbb{F}_n(t)(1-\mathbb{F}_n(t))}}}=0$ for $\mathbb{F}_n(t)\in\{0,1\}$.
The extreme value approximation (see \cite{Eicker} and \cite{J})
\begin{gather*}
\lim_{n\to \infty}\Pb_F\left(a_n \hat{T}_n -b_n\leq x\right) =\exp(-4\exp(-x)),\quad x\in\mathbb{R},
\end{gather*}
where $a_n$ and $b_n$ are as in (\ref{j1}), yields
the Eicker-Jaeschke confidence band $[\tilde{L}_n(t), \tilde{U}_n(t)]$ for $F(t)$ on the interval $t\in[X_{(1)}, X_{(n)})$ with the lower and upper bounds given by
\begin{gather*}
\tilde{L}_n(t)=\max\left\{0,\mathbb{F}_n(t)-a_n^{-1}({b_n+x_{\alpha}})\sqrt{{\mathbb{F}_n(t)(1-\mathbb{F}_n(t))}/{n}}\right\},\\
\tilde{U}_n(t)=\min\left\{1,\mathbb{F}_n(t)+a_n^{-1}({b_n+x_{\alpha}})\sqrt{{\mathbb{F}_n(t)(1-\mathbb{F}_n(t))}/{n}}\right\},
\end{gather*}
and $x_\alpha=-\log\left(-\log(1-\alpha)/4\right)$.

Numerical simulations show  that, even for moderate sample sizes,
when compared to the Kolmogorov-Smirnov confidence band, the CsCsHM confidence band is
of the same length ``in the middle'' and is shorter on the tails.
The new CsCsHM confidence band outperforms the Eicker-Jaeschke confidence band  ``in the middle'' and
does a similar job on the tails (see Figure~\ref{figure1} on page~\pageref{figure1}).

It is known that, as $n\to \infty$ (see \cite{Bickel} and \cite{Gnedenko}),
\begin{multline}\label{rate}
\sup_{-\infty<x<\infty}\left|\Pb_F\left(\sqrt{n}\sup_{-\infty<t<\infty}|\mathbb{F}_n(t)-F(t)| \leq x \right)\right.\\ -
\left.\Pb_F\left(\sup_{0<F(t)<1}|B(F(t)|\leq x \right)\right|
=O\left(n^{-1/2}\right),
\end{multline}
where $\{B(u),0\leq u\leq 1\}$ is a Brownian bridge. That is, under $H_0$, the CDF of the
two-sided Kolmogorov-Smirnov statistic converges to the Kolmogorov CDF $K(x)$, uniformly in $x\in\mathbb{R}$, at the rate of $O(n^{-1/2})$.
The results of numerical experiments suggest that the rate of convergence of the  CDF's of $T_n(q)$ and $T_n^+(q)$ to
their respective  limit CDF's  may be comparable to  that in (\ref{rate}). A theoretical justification of this claim is an open problem.

\begin{figure}[H]
\label{figure1}
\centering
\includegraphics[width=12cm,height=5.8cm]{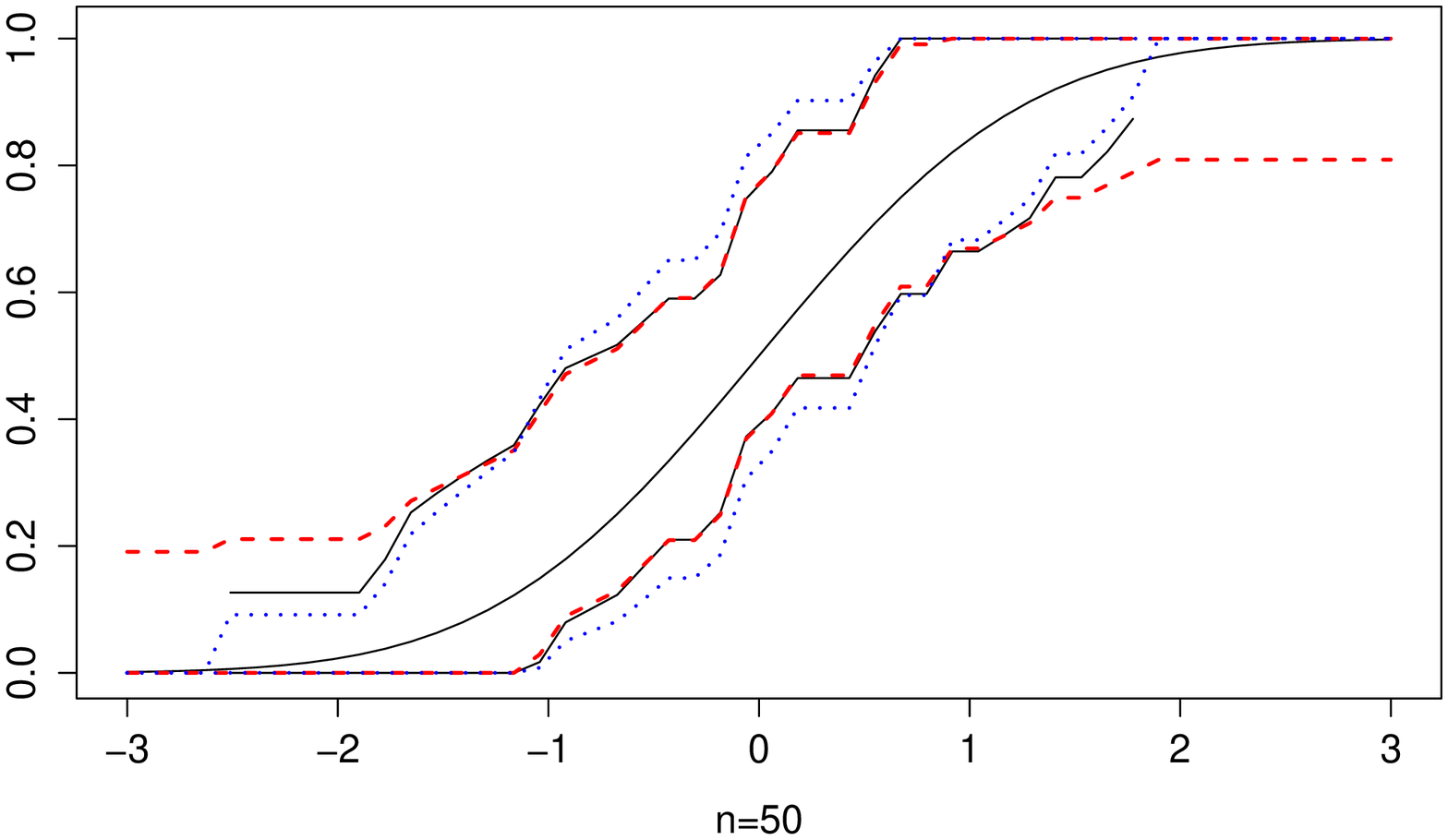}\\
\includegraphics[width=12cm,height=5.8cm]{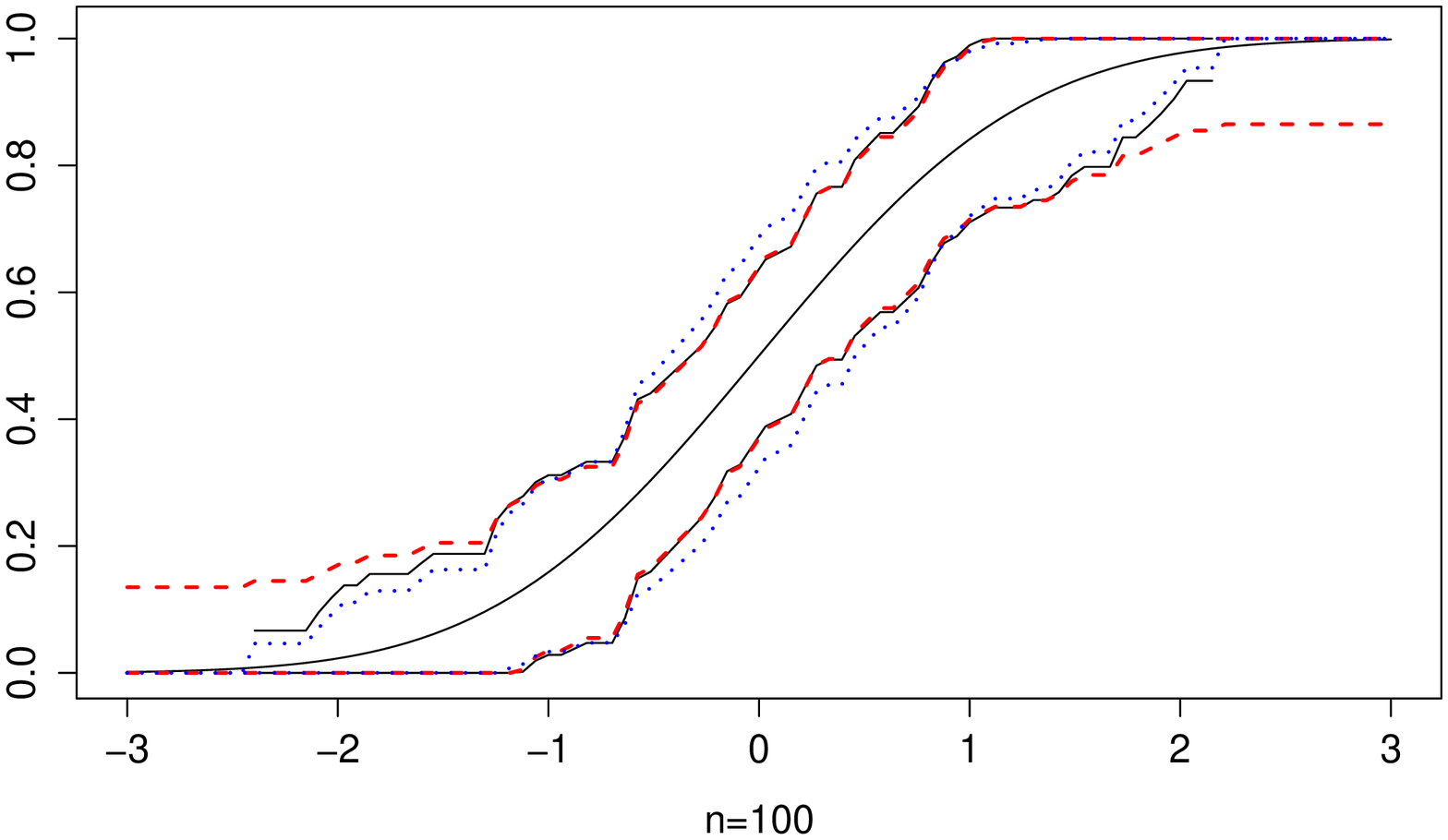}\\
\includegraphics[width=12cm,height=5.8cm]{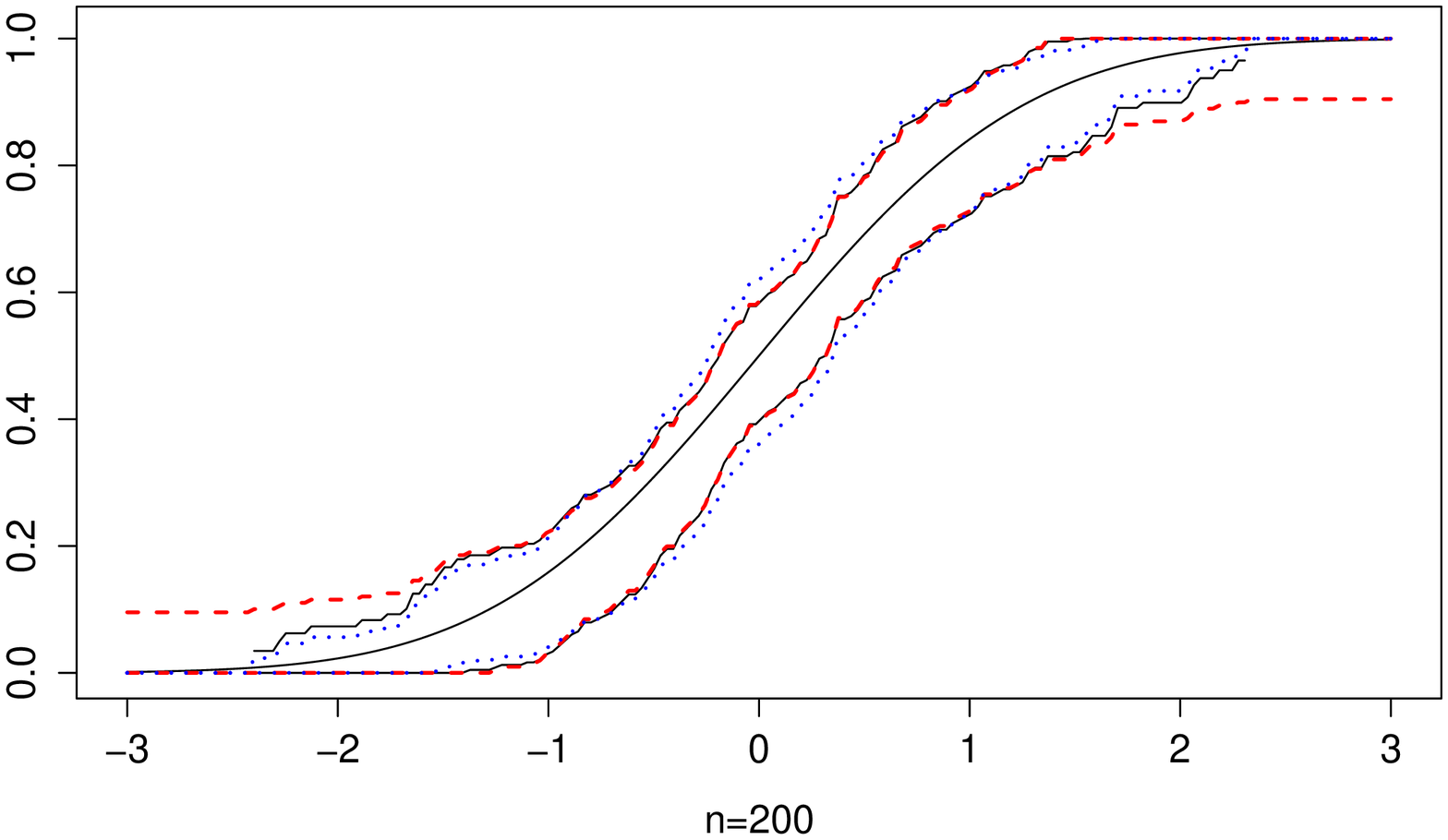}
\caption[]{Confidence bands for simulated data. The solid line is the true CDF.
The  solid lines above and below the middle line are a 95 percent Cs\"org\H{o}-Cs\"org\H{o}-Horv\'{a}th-Mason confidence band.
The red dashed lines are a 95 percent  Kolmogorov-Smirnov confidence band.
The blue dotted lines are a 95 percent Eicker-Jaeschke confidence band.}
\label{figure1}
\end{figure}

\section{Attainment of the Ingster optimal detection boundary}
\noindent
An important particular case of a goodness-of-fit testing problem is that of detecting sparse and weak heterogeneous mixtures.
The latter problem has been extensively studied after the publications of Ingster \cite{Ing97,Ing99}.
As in \cite{Ing97} (see also \cite{DJ, MR}), we first consider testing the null hypothesis
$$H_0: X_1,\ldots,X_n \,{\bd  {\mbox\small{{iid}}}\over \sim}\, N(0,1),$$
i.e., $F_0$ in (\ref{hnull}) is the standard normal CDF,
against a sequence of alternatives
$$H_{1,n}: X_1,\ldots,X_n\, {\bd  {\mbox\small{{iid}}}\over \sim} \, (1-\varepsilon_n)N(0,1)+\varepsilon_n N(\mu_n,1),$$
where $\varepsilon_n=n^{-\beta}$ for some \textit{sparsity index} $\beta\in(1/2,1)$ and $\mu_n=\sqrt{2r\log n}$ with $0<r<1$.
The parameters $\beta$ and $r$ are assumed unknown, and $n\to \infty$. The parameter $
\mu_n$ may be thought of as a \textit{signal strength}. In \cite{IPT} a similar mixture model emerged in connection
with a high-dimensional classification problem.

One well known property of a normal distribution says that if $\xi_1,\xi_2,\ldots$ is a sequence of iid standard normal random variables, then
$$P\left(\max_{1\leq i\leq n}|\xi_i|\geq \sqrt{2\log n}\right)\to 0,\quad n\to \infty.$$
This property explains the choice of the non-zero mean $\mu_n$ in the \textit{sparse} heterogeneous normal mixture specified by the alternative $H_{1,n}$: such a choice
makes the problem very hard but yet solvable.

In this section, we show
that if the parameter $r$ exceeds the \textit{detection boundary} $\rho(\beta)$ obtained by Ingster (see Section 2.6 of \cite{Ing97}; see also Section 1.1 of \cite{DJ}),
which is defined by
\begin{gather*}
\rho(\beta)=\left\{\begin{array}{ll} \beta-1/2,\quad\quad \;\,\,\quad 1/2<\beta< 3/4,\\ \vspace{-0.3cm}\\
(1-\sqrt{1-\beta})^2,\quad 3/4\leq\beta<1,
\end{array}\right.
\end{gather*}
then the test procedure based on  $T^+_n(q)$ distinguishes between
$H_0$ and $H_{1,n}$ (see Theorem 4.1). Since $T^+_n(q)$ does not require the knowledge of $\beta$ and $r$, following Donoho and Jin \cite{DJ},
we will call such a test procedure  \textit{optimally adaptive.}

Another model of interest, which was found to be useful in various classification problems (see, for example, \cite{Pavlenko12} and
\cite{Tillander13}), has the form:
\begin{eqnarray*}
H^{\prime}_0 &:& X_1,\ldots,X_n \,{\bd  {\mbox\small{{iid}}}\over \sim}\, \chi^2_\nu(0),\\
H_{1,n}^{\prime}&:& X_1,\ldots,X_n \,{\bd  {\mbox\small{{iid}}}\over \sim}\, (1-\varepsilon_n)\chi^2_\nu(0)+\varepsilon_n\chi^2_\nu(\delta_n),
\end{eqnarray*}
where  $\chi^2_\nu(\delta)$ denotes the noncentral chi-square distribution with $\nu$ degrees of freedom and
 noncentrality parameter $\delta$, $\varepsilon_n=n^{-\beta}$ for some $\beta\in(1/2,1)$, and $\delta_n=2r\log n$ for some $0<r<1.$
For $\nu=2$ this model is connected to the problem of detecting covert communications (see Section 1.7 in \cite{DJ}).

In order to apply the previously developed theory to the problem of testing $H_0$ ($H^{\prime}_0$) versus $H_{1,n}$ ($H_{1,n}^{\prime}$),
we need to transform the initial observations.
Namely, let $Y_i=1-\Phi(X_i)$ and let ${\cal G}(u)$ denote a common CDF of  the $Y_i$'s taking values in $[0,1]$.
Then the problem of testing $H_0$ versus $H_{1,n}$ transforms to testing
$${\cal H}_0: {\cal G}(u)= F_0(u),\quad \mbox{the uniform $U(0,1)$ CDF}$$
against a sequence of upper-tailed alternatives
$${\cal H}_{1,n}: {\cal G}(u)=F_0(u)+\varepsilon_n\left((1-u)-\Phi\left(\Phi^{-1}(1-u)-\mu_n \right) \right)>F_0(u).$$
The test statistic takes the form
\begin{gather*}
{T}^+_n(q)=\sup_{0<u<1}\frac{\sqrt{n}(\mathbb{G}_n(u)-u)}{q(u)},
\end{gather*}
where $\mathbb{G}_n(u)=n^{-1}\sum_{i=1}^n \mathbb{I}(Y_i\leq u)$ is the EDF based on the transforms variables $Y_i$'s.

In a chi-square mixture model, let $S_i=1-H_{\nu,0}(X_i)$, where $H_{\nu,\delta}$ is the CDF of a $\chi^2_{\nu}(\delta)$ distribution,
and let ${\cal H}(u)$ denote a common CDF of the $S_i$'s. Then
the problem of testing $H^{\prime}_0$ versus $H^{\prime}_{1,n}$ transforms to testing
$${\cal H}^{\prime}_0: {\cal H}(u)= F_{0}(u),\quad \mbox{the uniform $U(0,1)$ CDF}$$
against a sequence of upper-tailed alternatives
$${\cal H}^{\prime}_{1,n}: {\cal H}(u)=F_0(u)+\varepsilon_n\left((1-u)-H_{\nu,\delta_n}\left(H_{\nu,0}^{-1}(1-u)\right) \right)>F_0(u).$$
The test statistic becomes
\begin{gather*}
{T}^+_n(q)=\sup_{0<u<1}\frac{\sqrt{n}(\mathbb{H}_n(u)-u)}{q(u)},
\end{gather*}
where $\mathbb{H}_n(u)=n^{-1}\sum_{i=1}^n \mathbb{I}(S_i\leq u)$ is the EDF based on the $S_i$'s.

In both (normal and chi-square) settings, where we have data which are ``sparsely non-null'', our statistic $T_n^+(q)$ will be shown to be optimally adaptive
(see Theorems 4.1 and 4.2 below).

In connection with testing ${\cal H}_0$ versus ${\cal H}_{1,n}$ using the test statistic $T_n^+(q)$ with an EFKP upper-class function
$q$ we have the following result.

\medskip

\noindent \textbf{Theorem 4.1.}
\textit{For a function $q$ as in {\rm (\ref{qq})}, consider the test of asymptotic level $\alpha$ that rejects ${\cal H}_0$ when
$$T^+_n(q)\geq t^+_\alpha(q),$$
where the critical value $t^+_\alpha(q)$ is chosen to have $\Pb\left(\sup_{0<u<1}B(u)/q(u)\geq t^+_\alpha(q)\right)=\alpha.$
For every alternative ${\cal H}_{1,n}$ with $r$ exceeding the detection boundary $\rho(\beta)$,
the asymptotic level $\alpha$ test based on $T_n^+(q)$ has a full power, that is,}
$$\Pb_{{\cal H}_{1,n}}(T_n^+(q)\geq t^+_\alpha(q))\to 1,\quad n\to \infty.$$

\medskip

In connection with testing ${\cal H}^{\prime}_0$ versus ${\cal H}_{1,n}^{\prime}$ using the test statistic $T_n^+(q)$ with an EFKP upper-class function
$q$, we have the following result.

\medskip

\noindent \textbf{Theorem 4.2.}
\textit{For a function $q$ as in {\rm (\ref{qq})}, consider the test of asymptotic level $\alpha$ that rejects ${\cal H}^{\prime}_0$ when
$$T^+_n(q)\geq t^+_\alpha(q),$$
where the critical value $t^+_\alpha(q)$ is as in Theorem 4.1.
For every alternative ${\cal H}_{1,n}^{\prime}$ with $r$ exceeding the detection boundary $\rho(\beta)$,
the asymptotic level $\alpha$ test based on $T_n^+(q)$ has a full power, that is,}
$$\Pb_{{\cal H}^{\prime}_{1,n}}(T_n^+(q)\geq t^+_\alpha(q))\to 1,\quad n\to \infty.$$

\medskip

Theorems 4.1 and 4.2 say that if  $r>\rho(\beta)$, then asymptotically our test procedure based on $T_n^+(q)$ distinguishes between ${\cal H}_0$ and ${\cal H}_{1,n}$,
as well as between ${\cal H}^{\prime}_0$ and ${\cal H}_{1,n}^{\prime}$.
The proofs of both results are given in Section 7.

\medskip

\noindent \textit{Remark 4.1.}
{In view of Proposition 3.1, results similar to Theorems 4.1 and 4.2  hold true for the
whole class of statistics $T^+_n(q,I)$ indexed by a subinterval $I=(a,b)\subseteq(0,1)$,
in which case the critical region takes the form
$$T^+_n(q,I)\geq t^+_{\alpha}(q,I),$$
where  $t^+_{\alpha}(q,I)$  is determined by $\Pb(\sup_{a<u<b}B(u)/q(u)\geq  t^+_{\alpha}(q,I))=\alpha.$
In particular, this observation applies to the statistic $T^+_{n}(q,(0,\alpha_0))$ as in {\rm(\ref{nasha})}.}

\medskip

\noindent \textit{Remark 4.2.}
{Theorems 4.1 and 4.2 remain valid when the weight function $q$ defined  as in {\rm (\ref{qq})} is replaced by the Chibisov--O'Reilly function
$$q(u)=\left(u(1-u)\right)^{1/2}\left(\log\log(1/(u(1-u)))\right)^{1/2+\sigma},\quad \sigma>0.$$
The statement follows from the proof of Theorem 4.1 in Section 7, with $W_n(s)$  in {\rm(\ref{Wn})} replaced by
$$W_n(s)= {V_n(s)}/{\left(\log\log({1}/(p_{n,s}(1-p_{n,s})))\right)^{1/2+\sigma}},$$
and the subsequent derivations adjusted accordingly.
}

\section{Further remarks and comments}

\noindent

In recent years, some of the researchers who have begun to focus on detecting sparse heterogeneous mixtures,
have strongly advocated the use of the higher criticism statistic ${\rm HC}_n^+$ as in (\ref{DJT}) and its modifications
(see, for example, \cite{Cai, DJ} and Section 34.7 in \cite{DasG}).
In all these studies, attempts to numerically justify  theoretical properties of ${\rm HC}_n^+$
have resulted in sample sizes like $n=10^6$ and greater.
We wish to note, however, that this is due to the fact that, under $H_0$, the statistics ${\rm HC}_n$, ${\rm HC}^+_n$, and ${\rm HC}_n^*$ tend to $\infty$ in probability
(see \cite{Eicker} and \cite{J}), as well as almost surely (see Chapter 16 in \cite{Shorak} and references therein).
This property of the higher criticism statistic complicates its use in practice.

As noticed in the literature, one disadvantage of using ${\rm HC}^+_n$ is that
one has no clear recipe for the choice of its critical value.
Indeed, the test based on $\mbox{HC}^{+}_n$ prescribes to reject $H_0$ in favour of $H_{1,n}$  when
$$\mbox{HC}^{+}_n>h(n,\alpha_n),$$
where $h(n,\alpha_n)=\sqrt{2\log\log n}(1+o(1))$ and the level $\alpha_n\to 0$ slowly enough.
The problem of determining a critical value is unavoidable because, as mentioned just above,  ${\rm HC}^+_n$ tends almost surely to $\infty$  under $H_0$.
Unfortunately, this is not the only problem with applying the higher criticism in practice (see Section 2.3 for details).

We must  emphasize that our test procedure based on $T_n^+(q)$ is of a different kind.
This is due to the proper choice of the weight function $q(u)$, such as defined in (\ref{qq}),
which for all large enough $n$ makes the value of the sup-functional $\sup_{0<u<1} {\sqrt{d} (\mathbb{U}_{d}(u)-u)}/{q(u)}$
finite almost surely (see (\ref{fv}) and (\ref{mc2})).
Specifically, our test procedure rejects the null hypothesis at asymptotic level $\alpha$ when
$$T_n^+(q)>t^+_\alpha(q),$$
where the critical value  $t^+_\alpha(q)$ satisfies $\Pb\left(\sup_{0<u<1}B(u)/q(u)\geq t^+_\alpha(q)\right)=\alpha,$ and the distribution of
$\sup_{0<u<1}B(u)/q(u)$ is tabulated in Table~\ref{table1} on page~\pageref{table1}.

By proposing to use $T_n^+(q)$ or, more generally, a class of the CsCsHM-type test statistics $T_n^+(q,I)$, $I=(a,b)\subseteq(0,1)$,
instead of the higher criticism statistic, we are aiming at two goals. First, by using such a class, we obtain test
statistics that are {\it sensitive} to the choice of $\alpha_{0}$. This gives a correct implementation of the initial idea of Donoho and Jin \cite{DJ}
 that originates from the Tukey's concept of second-level significance testing in a multiple hypothesis testing setup.  Second, with the
class of CsCsHM-type test statistics available, we have an analytical solution to the ``end-points problem'', which eliminates the need of
 arbitrarily truncating the interval $(0,\alpha_0)$,
over which the supremum in (\ref{DJ}) is taken. Moreover,
in various signal detection problems involving unknown parameters, the tests based on $T_n^+(q,I)$, $I=(a,b)\subseteq(0,1)$, are found to be optimally adaptive.

The main results of this paper, Propositions 2.1 and 3.1 and Theorems 4.1 and 4.2, show that,
in the sup-norm scenario, when normalizing the empirical process $\sqrt{n}|\mathbb{U}_n(u)-u|$ by an EFKP upper-class function
$q(u)$, we do exactly the right job. Our conjecture therefore is that, in case of some other non-Gaussian heterogeneous mixtures,
including those studied in Section 5 of \cite{DJ},
similar results on the distinguishability of the null and alternative hypotheses by means of the test statistics $T^+_n(q,I)$, $I=(a,b)\subseteq(0,1)$, remain valid.

\section{Tabulation of cumulative distribution functions}

\noindent
It follows from the analysis of the previous sections that the CsCsHM test statistics with EFKP weight functions have a number of attractive features
that could be very useful in practical applications.
In particular, the limit distributions of these statistics under the null hypothesis are easily tabulated.

In this section, we assume that $$q(u)=\sqrt{u(1-u)\log\log(1/u(1-u))},\quad 0<u<1.$$
The distribution of the random variable $\sup_{0<u<1}{|B(u)|}/{q(u)}$ has been tabulated
in B. Eastwood and V. Eastwood \cite{EE} and, using a somewhat different approach, in Orasch and Pouliot \cite{OP}.
In this section, the tabulation of the distribution of  $\sup_{0<u<1}{B(u)}/{q(u)}$ follows the approach of
\cite{OP}. Namely, we shall use the following algorithm.

\begin{enumerate}
\item[{1.}] Choose a large positive integer $n$. Generate $n$ independent normal
 $N(0,1)$ random variables $X_1,\ldots,X_n$.

\item[{2.}] Choose a large positive integer $M$. Repeat step 1 $M$ times, and for
$m=1,\ldots, M$, let $X_1^{(m)},\ldots,X_n^{(m)}$ denote the  data obtained on the
$m$th iteration.

\item[{3.}] For each $m=1,\ldots, M$, calculate the partial sums
$S^{(m)}_k=\sum_{i=1}^k X_i^{(m)},$ $k=1,\ldots,n.$

\item[{4.}] For each $m=1,\ldots, M$, find the value of
$$T_n^{(m)}=\max_{1\leq k\leq n-1}\frac{S_k^{(m)}-(k/n)S_n^{(m)}}{q(k/n)n^{1/2}}.$$

\item[{5.}] For $x\in\mathbb{R}$, use the function
$G_{n,M}(x)=\frac{1}{M}\sum_{m=1}^M \mathbb{I}\left( T_n^{(m)}\leq x \right)$
to approximate the limit CDF $G(x)=\Pb\left(\sup\limits_{0<u<1}{B(u)}/{q(u)}\leq x\right)$.

\end{enumerate}

\medskip
\noindent

\textit{Remark 6.1.}
{In step 1 of the algorithm, $X_1,\ldots,X_n$ could be a random sample from any distribution with
$\Eb(X_1^2)<\infty$.
Then, one would have to standardize $T_n^{(m)}$ in step 4 by putting $\sigma=\sqrt{\Var (X_1)}$ into the denominator.
The finiteness of $\Eb(X_1^2)$ makes it possible to apply part (ii) of Theorem 2.1.1 in \cite{CsHorvath},
which, together with the Glivenko-Cantelli theorem,
guaranties the closeness of $G_{n,M}(x)$ and $G(x)$ for all large enough  $n$ and $M$.}

\medskip

Step 5 of the algorithm is based on the following result of \cite{CsHorvath}
that is similar to the result in Fact 1.
 Namely, let $X_1,\ldots,X_n$ be a random sample from a distribution with
 $\Eb(X_1^2)<\infty$. Consider the normalized tied-down partial sums process
 \begin{gather*}
Z_n(u)= \left\{\begin{array}{ll} \left(S_{[(n+1)u]}-[(n+1)u]S_n/n \right)/(n^{1/2}\sigma),\quad\; \;\; 0\leq u<1,\\ \vspace{-0.3cm}\\
0 ,\qquad \qquad\qquad \qquad\qquad\qquad \qquad \qquad\quad \,\, \,u=1,
\end{array}\right.
 \end{gather*}
 where $\sigma^2=\Var(X_1)$.
Then, by part (ii) of Theorem 2.1.1 in \cite{CsHorvath} (see also part (b) of Corollary 2.1 in \cite{CsHor}),
\begin{gather}\label{cor212}
\sup_{0<u<1}Z_n(u)/q(u)
\,\,{\bd  {\mbox\small{{\cal D}}}\over \rightarrow}\,\, \sup_{0<u<1} B(u)/q(u),
\quad n\to \infty,
\end{gather}
and
\begin{gather*}
\sup_{0<u<1}|Z_n(u)|/q(u)
\,\,{\bd  {\mbox\small{{\cal D}}}\over \rightarrow}\,\, \sup_{0<u<1} |B(u)|/q(u),
\quad n\to \infty,
\end{gather*}
if and only if $q$ is an EFKP upper-class function. An obvious modification of this result,
with the supremum over an arbitrary interval $(a,b)$, $0\leq a<b\leq 1$, also holds true.

Now, for the random sample $X_1,\ldots, X_n$  generated in step 1, define the function
$$G_n(x)=\Pb\left( \max_{1\leq k\leq n-1}\frac{S_k-(k/n)S_n}{q(k/n)n^{1/2}} \leq x \right),\quad x\in\mathbb{R},$$
where $S_k=\sum_{i=1}^k {X_i},$ and observe that by the triangle inequality, for every $x\in\mathbb{R},$
\begin{gather}\label{eq1}
|G(x)-G_{n,M}(x)|\leq |G(x)-G_n(x)|+|G_{n}(x)-G_{n,M}(x)|.
\end{gather}
It follows from the Glivenko-Cantelli theorem that for all $n\geq 2$,
\begin{gather}\label{eq2}
\lim_{M\to \infty}\sup_{x\in\mathbb{R}}|G_n(x)-G_{n,M}(x)| \,\,{\bd  {\mbox\small{{\rm{a.s.}}}}\over =}\,\,  0.
\end{gather}
Next, by means of (\ref{cor212}), as $n\to \infty$
\begin{gather*}
\max_{1\leq k\leq n-1}\frac{S_k-(k/n)S_n}{q(k/n)n^{1/2}} \,\,{\bd  {\mbox\small{{\cal D}}}\over \rightarrow}\,\, \sup_{0<u<1} B(u)/q(u).
\end{gather*}
From this, using the fact that $G(x)$ is a continuous function,
\begin{gather}\label{eq3}
\lim_{n\to \infty}\sup_{x\in\mathbb{R}}|G(x)-G_{n}(x)| = 0.
\end{gather}
It now follows from (\ref{eq1})--(\ref{eq3}) that  
$$\lim_{n\to\infty\atop M\to \infty}\sup_{x\in\mathbb{R}}| G(x)-G_{n,M}(x)| \,\,{\bd  {\mbox\small{{\rm{ a.s.}}}}\over =}\,\, 0.$$

Table~\ref{table1}  on page~\pageref{table1}  contains percentage points of the distribution function $G(x)$. For tabulating $G(x)$
the above algorithm with $n=M=50,000$ has been applied. From Table~\ref{table1}, the upper $1\%$, $5\%$, and $10\%$ percentage points are
5.16, 4.14, and 3.62, respectively. Note that Table~\ref{table1} supports the analytical finding of Proposition 3.1, according to which the tails
have been tamed. The modification of this algorithm to the case of test statistic $T_n^+(q,I)$ depending on a subinterval $I\subseteq(0,1)$
 is obvious.

\begin{table*}
\caption{The limit distribution of $\sup\limits_{0<u<1}\dfrac{\sqrt{n}(\mathbb{U}_n(u)-u)}{\sqrt{u(1-u)\log\log(1/u(1-u))}}$.
}
\label{table1}
\medskip
{\centerline
{\begin{tabular}{cccccc}
 \hline \\ [-7pt]
$x$ \; & $G(x)\;\;\;\;\;\;\;$ & $x$ \;& $G(x)\;\;\;\;\;\;\;$ & $x$\; & $G(x)$ \\ [5pt]
\hline \\ [-5pt]
0.74 \;& 0.01 \;\;\;\;\;\;\; & 1.81 \;& 0.34 \;\;\;\;\;\;\; & 2.57 \;& 0.67\\
0.87 \;& 0.02 \;\;\;\;\;\;\; & 1.83 \;& 0.35 \;\;\;\;\;\;\; & 2.60 \;& 0.68\\
0.95 \;& 0.03 \;\;\;\;\;\;\; & 1.85 \;& 0.36 \;\;\;\;\;\;\; & 2.63 \;& 0.69\\
1.02 \;& 0.04 \;\;\;\;\;\;\; & 1.87 \;& 0.37 \;\;\;\;\;\;\; & 2.66 \;& 0.70\\
1.07 \;& 0.05 \;\;\;\;\;\;\; & 1.89 \;& 0.38 \;\;\;\;\;\;\; & 2.69 \;& 0.71\\
1.11 \;& 0.06 \;\;\;\;\;\;\; & 1.91 \;& 0.39 \;\;\;\;\;\;\; & 2.72 \;& 0.72\\
1.16 \;& 0.07 \;\;\;\;\;\;\; & 1.93 \;& 0.40 \;\;\;\;\;\;\; & 2.76 \;& 0.73\\
1.19 \;& 0.08 \;\;\;\;\;\;\; & 1.95 \;& 0.41 \;\;\;\;\;\;\; & 2.79 \;& 0.74\\
1.23 \;& 0.09 \;\;\;\;\;\;\; & 1.97 \;& 0.42 \;\;\;\;\;\;\; & 2.83 \;& 0.75\\
1.26 \;& 0.10 \;\;\;\;\;\;\; & 1.99 \;& 0.43 \;\;\;\;\;\;\; & 2.87 \;& 0.76\\
1.29 \;& 0.11 \;\;\;\;\;\;\; & 2.01 \;& 0.44 \;\;\;\;\;\;\; & 2.91 \;& 0.77\\
1.32 \;& 0.12 \;\;\;\;\;\;\; & 2.03 \;& 0.45 \;\;\;\;\;\;\; & 2.95 \;& 0.78\\
1.35 \;& 0.13 \;\;\;\;\;\;\; & 2.05 \;& 0.46 \;\;\;\;\;\;\; & 2.99 \;& 0.79\\
1.37 \;& 0.14 \;\;\;\;\;\;\; & 2.07 \;& 0.47 \;\;\;\;\;\;\; & 3.03 \;& 0.80\\
1.40 \;& 0.15 \;\;\;\;\;\;\; & 2.09 \;& 0.48 \;\;\;\;\;\;\; & 3.08 \;& 0.81\\
1.42 \;& 0.16 \;\;\;\;\;\;\; & 2.12 \;& 0.49 \;\;\;\;\;\;\; & 3.13 \;& 0.82\\
1.45 \;& 0.17 \;\;\;\;\;\;\; & 2.14 \;& 0.50 \;\;\;\;\;\;\; & 3.18 \;& 0.83\\
1.47 \;& 0.18 \;\;\;\;\;\;\; & 2.16 \;& 0.51 \;\;\;\;\;\;\; & 3.23 \;& 0.84\\
1.49 \;& 0.19 \;\;\;\;\;\;\; & 2.18 \;& 0.52 \;\;\;\;\;\;\; & 3.29 \;& 0.85\\
1.51 \;& 0.20 \;\;\;\;\;\;\; & 2.20 \;& 0.53 \;\;\;\;\;\;\; & 3.35 \;& 0.86\\
1.54 \;& 0.21 \;\;\;\;\;\;\; & 2.22 \;& 0.54 \;\;\;\;\;\;\; & 3.42 \;& 0.87\\
1.56 \;& 0.22 \;\;\;\;\;\;\; & 2.25 \;& 0.55 \;\;\;\;\;\;\; & 3.48 \;& 0.88\\
1.58 \;& 0.23 \;\;\;\;\;\;\; & 2.27 \;& 0.56 \;\;\;\;\;\;\; & 3.55 \;& 0.89\\
1.60 \;& 0.24 \;\;\;\;\;\;\; & 2.30 \;& 0.57 \;\;\;\;\;\;\; & 3.62 \;& 0.90\\
1.63 \;& 0.25 \;\;\;\;\;\;\; & 2.32 \;& 0.58 \;\;\;\;\;\;\; & 3.70 \;& 0.91\\
1.65 \;& 0.26 \;\;\;\;\;\;\; & 2.35 \;& 0.59 \;\;\;\;\;\;\; & 3.79 \;& 0.92\\
1.67 \;& 0.27 \;\;\;\;\;\;\; & 2.37 \;& 0.60 \;\;\;\;\;\;\; & 3.89 \;& 0.93\\
1.69 \;& 0.28 \;\;\;\;\;\;\; & 2.40 \;& 0.61 \;\;\;\;\;\;\; & 4.00 \;& 0.94\\
1.71 \;& 0.29 \;\;\;\;\;\;\; & 2.43 \;& 0.62 \;\;\;\;\;\;\; & 4.14 \;& 0.95\\
1.73 \;& 0.30 \;\;\;\;\;\;\; & 2.46 \;& 0.63 \;\;\;\;\;\;\; & 4.30 \;& 0.96\\
1.75 \;& 0.31 \;\;\;\;\;\;\; & 2.49 \;& 0.64 \;\;\;\;\;\;\; & 4.48 \;& 0.97\\
1.77 \;& 0.32 \;\;\;\;\;\;\; & 2.51 \;& 0.65 \;\;\;\;\;\;\; & 4.73 \;& 0.98\\
1.79 \;& 0.33 \;\;\;\;\;\;\; & 2.54 \;& 0.66 \;\;\;\;\;\;\; & 5.16 \;& 0.99
\end{tabular}}}
\end{table*}

\section{Proofs} This section contains the proofs of Proposition 2.1 and Theorems 4.1 and 4.2 from the previous sections.

\medskip
\noindent {\textit{Proof of Proposition 2.1.}} The proof of this result goes long the lines of Section 4.4 in \cite{CsCsHM}.
We show that statement (\ref{j2}) follows from
the result of Darling and Erd\H{o}s \cite{Darling} cited below.

\medskip

\noindent \textbf{Fact 3 (Theorem 1.9.1 (Darling and Erd\H{o}s \cite{Darling}) in \cite{CsR}).} \textit{Let
 $\{U(t),-\infty<t<\infty\}$ be the Ornstein--Uhlenbeck process, and let
$$a(y,T)=(y+2\log T+(1/2)\log\log T-(1/2)\log \pi)(2\log T)^{-1/2},\quad y\in\mathbb{R},\quad T>1.$$
Then}
\begin{gather*}
\lim_{T\to \infty}\Pb\left( \sup_{0\leq t\leq T} U(t)\leq a(y,T) \right)=\exp(-\exp(-y)),\\
\lim_{T\to \infty}\Pb\left( \sup_{0\leq t\leq T} |U(t)|\leq a(y,T) \right)=\exp(-2\exp(-y)),
\end{gather*}

Before using Fact 3, observe that
\begin{gather*}
\{U(t),-\infty<t<\infty\} \,\,{\bd  {\mbox\small{{\cal D}}}\over =}\,\, \left\{(1+e^{2t})e^{-t} B\left(\frac{e^{2t}}{1+e^{2t}}\right),-\infty<t<\infty   \right\},
\end{gather*}
where $\{B(u),0\leq u\leq 1\}$ is a Brownian bridge. Then, using the stationarity of the Ornstein--Uhlenbeck process $U(t)$, we have
for any decreasing sequence of numbers $\varepsilon_n\to 0$, any number $\alpha_0\in(0,1)$, and any $y\in\mathbb{R}$,
\begin{gather*}
\lim_{n\to \infty}\Pb\left\{\sup_{\varepsilon_n<u<\alpha_0}\frac{B(u)}{\sqrt{u(1-u)}}\leq a\left(y,\frac12\log\frac{\alpha_0(1-\varepsilon_n)}{\varepsilon_n(1-\alpha_0)}\right)\right\}\\
=\lim_{n\to \infty}\Pb\left\{\sup_{\frac12\log\frac{\varepsilon_n}{1-\varepsilon_n}<t<\frac12\log\frac{\alpha_0}{1-\alpha_0}}U(t)\leq a\left(y,\frac12\log\frac{\alpha_0(1-\varepsilon_n)}{\varepsilon_n(1-\alpha_0)}\right)\right\}\\
=\lim_{n\to \infty}\Pb\left\{\sup_{0<t<\frac 12\log\frac{\alpha_0(1-\varepsilon_n)}{\varepsilon_n(1-\alpha_0)}}U(t)\leq a\left(y,\frac12\log\frac{\alpha_0(1-\varepsilon_n)}{\varepsilon_n(1-\alpha_0)}\right)\right\}.
\end{gather*}
Then, by Fact 3, cf. Lemma 4.4.1 in \cite{CsCsHM},
\begin{gather}\label{j3}
\lim_{n\to \infty}\Pb\left\{\sup_{\varepsilon_n<u<\alpha_0}\frac{B(u)}{\sqrt{u(1-u)}}\leq a\left(y,\frac12\log\frac{\alpha_0(1-\varepsilon_n)}{\varepsilon_n(1-\alpha_0)}\right)\right\}
=\exp(-\exp(-y)).
\end{gather}
Next, choose
$$\varepsilon_n=(\log n)^3/n$$
and observe that for any $0<\alpha_0<1$ and any $y\in\mathbb{R}$, 
as $n\to \infty$
\begin{gather*}
a\left(y,\frac12\log\frac{\alpha_0(1-\varepsilon_n)}{\varepsilon_n(1-\alpha_0)}\right)= a\left(y,(\log n)/2)(1+o(1))\right)=
a\left(y, (\log n)/2\right)(1+o(1)),
\end{gather*}
where
\begin{gather*}
a\left(y,(\log n)/2\right)=\frac{y+2\log \left((\log n)/2\right)+(1/2)\log\log\left(( \log n)/2\right)-(1/2)\log \pi }{\sqrt{2\log \left((\log n)/2\right)}}\\
=\frac{y+2\log \log n +(1/2)\log \log \log n-(1/2)\log \pi+o(1)}{\sqrt{(2\log \log n) (1+o(1))}}\\
= \frac{y-\log 2+b_n+o(1)}{a_n(1+o(1))},
\end{gather*}
with $a_n$ and $b_n$ as in (\ref{j1}). From this,
introducing in (\ref{j3}) a new variable $x$  by the formula $x=y-\log 2$, we can write
\begin{gather}\label{j4}
\lim_{n\to \infty}\Pb\left\{a_n\sup_{\varepsilon_n<u<\alpha_0}\frac{B(u)}{\sqrt{u(1-u)}}-b_n\leq x\right\}
=\exp\left(-\frac12\exp(-x)\right).
\end{gather}

Now define
\begin{eqnarray*}
V_n&=&a_n\sup_{0<u<\alpha_0}\frac{\sqrt{n}(\mathbb{U}_n(u)-u)}{\sqrt{u(1-u)}}-b_n,\\
V^{(1)}_n&=&a_n\sup_{0<u\leq\varepsilon_n}\frac{\sqrt{n}(\mathbb{U}_n(u)-u)}{\sqrt{u(1-u)}}-b_n,\\
V^{(2)}_n&=&a_n\sup_{\varepsilon_n<u<\alpha_0}\frac{\sqrt{n}(\mathbb{U}_n(u)-u)}{\sqrt{u(1-u)}}-b_n.
\end{eqnarray*}
Then
$$V_n=\max(V_n^{(1)}, V_n^{(2)}),$$
where by (4.4.22) in \cite{CsCsHM}
\begin{gather}\label{star}
V_n^{(1)} \,\,{\bd  {\mbox\small{P}}\over \rightarrow}\,\, -\infty.
\end{gather}
Hence the limit distribution of $V_n$ coincides with that of $V_n^{(2)}$.

Now observe that the probability space on which the random variable $U_i$'s are defined can be extended in such a way
that on the new (extended) probability space we can construct a sequence of Brownian bridges
$\{B_n(u),0\leq u\leq 1\}$, $n=1,2,\ldots,$ such that for any $0<\alpha_0<1$ and any $0<\nu<1/4$, as $n\to \infty$
(see formula (4.2.5b) in Corollary 4.2.2 of \cite{CsCsHM})
\begin{gather}\label{bd}
n^{\nu}\sup_{0<u<\alpha_0}\frac{{\left|\sqrt{n}(\mathbb{U}_n(u)-u)-\bar{B}_n(u)\right|}}{{(u(1-u))^{1/2-\nu}}}=O_P(1),
\end{gather}
where
\begin{gather*}
\bar{B}_n(u)=\left\{\begin{array}{ll}
B_n(u),\quad 1/n\leq u\leq 1-1/n,
\\ \vspace{-0.3cm}\\
0,\quad \quad \quad\mbox{elsewhere}.
\end{array}\right.
\end{gather*}
Next, for each $n$,
\begin{gather*}
\left|V_n^{(2)} -\left(a_n\sup_{\varepsilon_n< u<\alpha_0}\frac{B_n(u)}{\sqrt{u(1-u)}}-b_n\right)\right| \leq
a_n\sup_{\varepsilon_n<u<\alpha_0}\frac{\left|\sqrt{n}(\mathbb{U}_n(u)-u)-{B}_n(u)\right|}{\sqrt{u(1-u)}}\nonumber\\ \leq
\frac{2 a_n n^{\nu}}{(\log n)^{3\nu}}\sup_{1/n<u<\alpha_0}\frac{\left|\sqrt{n}(\mathbb{U}_n(u)-u)-{B}_n(u)\right|}{(u(1-u))^{1/2-\nu}},
\end{gather*}
which together with (\ref{bd}) yields as $n\to \infty$
\begin{gather*}
\left|V_n^{(2)} -\left(a_n\sup\limits_{\varepsilon_n\leq u<\alpha_0}\frac{B_n(u)}{\sqrt{u(1-u)}}-b_n\right)\right|=O_P\left(\frac{a_n}{(\log n)^{3\nu}}\right)=o_P(1). 
\end{gather*}
Therefore the limit distribution of $V_n^{(2)}$ (and hence of $V_n$) coincides with that of $K_n(B_n):= a_n\sup_{\varepsilon_n\leq u<\alpha_0}\frac{B_n(u)}{\sqrt{u(1-u)}}-b_n$.
Note that, for each $n$, $\{B_n(u), 0\leq u\leq 1\} \,\,{\bd  {\mbox\small{{\cal D}}}\over =}\,\,\{B(u),0\leq u\leq 1\}$ and hence
$K_n(B_n)\,\,{\bd  {\mbox\small{{\cal D}}}\over =}\,\, K_n(B)$. From this we get via relation (\ref{j4}) that
for every $x\in\mathbb{R},$
\begin{gather*}
\lim_{n\to \infty}\Pb\left(a_n \sup_{0<u<\alpha_0} \frac{\sqrt{n}(\mathbb{U}_n(u)-u)}{\sqrt{u(1-u)}}-b_n \leq x\right) =\exp\left(-\frac12\exp(-x)\right).
\end{gather*}
The proof is now complete.
\done
\medskip

\noindent {\textit{Proof of Theorem 4.1.}} The proof of this theorem goes along the lines of that of Theorem 1.2 in \cite{DJ}.
We need to show that
\begin{gather*}
\lim _{n\to \infty}\Pb_{{\cal H}_{1,n}}\left(T_n^+(q)\leq t^+_\alpha(q)  \right)= 0.
\end{gather*}

Let $0<s\leq 1 $ and introduce the following notation:
\begin{gather}
p_{n,s}=\Pb_{{\cal H}_0}\left(Y_i\leq \Phi\left(-\sqrt{2s\log n}\right)\right),\nonumber \\
p^{\prime}_{n,s}=\Pb_{{\cal H}_{1,n}}\left( Y_i\leq \Phi\left(-\sqrt{2s\log n}\right) \right),\nonumber \\
 N_n(s)=\#\left\{i:Y_i\leq \Phi\left(-\sqrt{2s\log n}\right)\right\},\nonumber\\
V_n(s)=\frac{N_{n}(s)-n p_{n,s}}{\sqrt{n p_{n,s}(1-p_{n,s})}},\nonumber\\
W_n(s)= \frac{N_{n}(s)-n p_{n,s}}{\sqrt{n p_{n,s}(1-p_{n,s})\log\log(1/(p_{n,s}(1-p_{n,s})))}}\nonumber\\
=\frac{V_n(s)}{\sqrt{\log\log({1}/(p_{n,s}(1-p_{n,s})))}}.\label{Wn}
\end{gather}
Since
\begin{gather*}
T_n^+(q)\geq \sup_{0<s\leq 1}W_n(s)\geq W_n(1)\, {\bd {\small{\cal D}}\over =}\, \frac{V_{n}(1)}{\sqrt{\log\log(1/(p_{n,1}(1-p_{n,1})))}}.
\end{gather*}
it follows that
\begin{gather*}
\Pb_{{\cal H}_{1,n}}\left(T_n^+(q)\leq t^+_\alpha(q)  \right)\leq \Pb_{{\cal H}_{1,n}}\left(W_n(1)\leq t^+_\alpha(q) \right)\\ =
\Pb_{{\cal H}_{1,n}}\left( V_n(1)\leq  t^+_\alpha(q) \sqrt{\log\log({1}/(p_{n,1}(1-p_{n,1})))}  \right)\\
= \Pb_{{\cal H}_{1,n}}\left( N_n(1)\leq n p_{n,1}+ t^+_\alpha(q)\sqrt{n p_{n,1}} \sqrt{\log\log({1}/(p_{n,1}(1-p_{n,1})))}  \right),
\end{gather*}
where, under ${\cal H}_{1,n}$, $N_n(1)$ is the sum iid Bernoulli random variables $Y_{i,n}$, $i=1,\ldots,n$, with parameter
$p^{\prime}_{n,1}$. 
Noting that
\begin{eqnarray*}
p_{n,s}&=&\Pb\left(N(0,1)\geq\sqrt{2s\log n}\right),\\
p^{\prime}_{n,s}&=&\Pb\left((1-\varepsilon_n)N(0,1)+\varepsilon_nN(\mu_n,1)\geq \sqrt{2s\log n} \right),
\end{eqnarray*}
and using the fact
$$\Pb\left(N(0,1)>x\right)\sim \frac{e^{-x^2/2}}{x\sqrt{2\pi}},\quad x\to \infty,$$
one can find that
\begin{eqnarray}\label{asympt}
p_{n,1}=O\left(n^{-1}\log^{-1/2}n\right),\quad
p^{\prime}_{n,1}=O\left( n^{-\beta-({1}-\sqrt{r})^2}\log^{-1/2}n \right).
\end{eqnarray}

\medskip
\textbf{Case 1.} Assume that either (a)  $3/4\leq\beta< 1$ and $r>\rho(\beta)=\left(1-\sqrt{1-\beta}\right)^2$ or
(b) $1/2<\beta< 3/4$ and $r\geq1/4$. Then, from the above
\begin{multline*}
\Pb_{{\cal H}_{1,n}}\left(T_n^+(q)\leq t^+_\alpha(q)  \right)\\ \leq
\Pb_{{\cal H}_{1,n}}\left( N_n(1)\leq n p_{n,1}+ t^+_\alpha(q)\sqrt{n p_{n,1}} \sqrt{\log\log({1}/(p_{n,1}(1-p_{n,1})))}  \right)\\
\leq \Pb\left(\sum_{i=1}^n(Y_{i,n}-p^{\prime}_{n,1}) \right.\leq - \left[ n p^{\prime}_{n,1}-n p_{n,1}-t^+_\alpha(q)\sqrt{n p_{n,1}}\times \right.
\\ \left.\left. \times\sqrt{\log\log({1}/(p_{n,1}(1-p_{n,1})))}  \right]\right),
\end{multline*}
where $Y_{1,n},\ldots,Y_{n,n}$ are iid Bernoulli random variables with parameter $p^{\prime}_{n,1}$ and (see (\ref{asympt}))
\begin{eqnarray*}
a_n&:=&n p^{\prime}_{n,1}-n p_{n,1}-t^+_\alpha(q)\sqrt{n p_{n,1}} \sqrt{\log\log({1}/(p_{n,1}(1-p_{n,1})))}\\
& =&O\left( n^{1-\beta-(1-\sqrt{r})^2}\log^{-1/2}n\right).
\end{eqnarray*}
Since, under our  assumptions on $\beta$ and $r$, the exponent $1-\beta-(1-\sqrt{r})^2$ is strictly positive, one has
$a_n\to \infty$ and also $a_n/\sqrt{n p^{\prime}_{n,1}}\to \infty$. Therefore, by using Chebyshev's inequality, we obtain
\begin{eqnarray*}
\lim_{n\to \infty}\Pb_{{\cal H}_{1,n}}\left(T_n^+(q)\leq t^+_\alpha(q)  \right)&\leq&
\lim_{n\to \infty}\Pb\left(\sum_{i=1}^n(Y_{i,n}-p^{\prime}_{n,1})\leq - a_n\right)\\&\leq&
\lim_{n\to \infty}\frac{n p^{\prime}_{n,1}}{a_n^2}=0.
\end{eqnarray*}

\medskip
\textbf{Case 2.} Assume that $1/2<\beta< 3/4$ and $\beta-1/2=\rho(\beta)<r<1/4$. Notice that in this case $\beta+r<1.$
Similar to Case 1, we can write
\begin{gather*}
\Pb_{{\cal H}_{1,n}}\left(T_n^+(q)\leq t^+_\alpha(q)  \right)\\ \leq \Pb_{{\cal H}_{1,n}}\left(\sup_{0<s\leq 1}W_n(s)\leq t^+_\alpha(q) \right)
\leq \Pb_{{\cal H}_{1,n}}\left(W_n(4r)\leq t^+_\alpha(q) \right)\\
=\Pb_{{\cal H}_{1,n}}\left( V_n(4r)\leq  t^+_\alpha(q) \sqrt{\log\log({1}/(p_{n,4r}(1-p_{n,4r})))}  \right)\\
 = \Pb_{{\cal H}_{1,n}}\left( N_n(4r)\leq n p_{n,4r}+ t^+_\alpha(q)\sqrt{n p_{n,4r}} \sqrt{\log\log({1}/(p_{n,4r}(1-p_{n,4r})))}  \right),
\end{gather*}
where, under ${\cal H}_{1,n}$, the random variable $N_n(4r)$ is the sum of iid Bernoulli random variables $Y^{\prime}_{1,n},\ldots,Y^{\prime}_{n,n}$ with parameter
$p^{\prime}_{n,4r}$. It is easily seen that
\begin{gather*}
p_{n,4r}=O\left( n^{-4r}\log^{-1/2}n\right),\quad
p^{\prime}_{n,4r}=p_{n,4r}+O\left( n^{-(\beta+r)}\log^{-1/2}n\right).
\end{gather*}
Therefore
\begin{gather*}
\lim_{n\to \infty}\Pb_{{\cal H}_{1,n}}\left(T_n^+(q)\leq t^+_\alpha(q)  \right)\leq
\lim_{n\to \infty}\Pb\left(\sum_{i=1}^n(Y^{\prime}_{i,n}-p^{\prime}_{n,4r})\leq - a^{\prime}_n\right),
\end{gather*}
where $a^{\prime}_n:=n(p^{\prime}_{n,4r}- p_{n,4r})-t^+_\alpha(q)\sqrt{np_{n,4r}}\sqrt{\log\log({1}/(p_{n,4r}(1-p_{n,4r})))}.$
Since
\begin{gather*}
a^{\prime}_n\asymp n^{1-(\beta+r)}\log^{-1/2}n, \quad n\to \infty,
\end{gather*}
it follows that
\begin{gather*}
\frac{a^{\prime}_n}{\sqrt{np^{\prime}_{n,4r}}} \asymp \left\{\begin{array}{ll} n^{r-(\beta-1/2)} \log^{-1/4}n ,\quad\quad\, \beta>3r,\\ \vspace{-0.3cm}\\
n^{1/2(1-(\beta+r))} \log^{-1/4}n ,\quad \beta\leq 3r,
\end{array}\right.
\end{gather*}
From this, under the above assumptions on the range of $r$ and $\beta$, we obtain
$$a^{\prime}_n\to \infty\quad\mbox{and}\quad {a^{\prime}_n}/{\sqrt{np^{\prime}_{n,4r}}}\to \infty.$$
By Chebyshev's inequality, this implies
\begin{eqnarray*}
\lim_{n\to \infty}\Pb_{{\cal H}_{1,n}}\left(T_n^+(q)\leq t^+_\alpha(q)  \right)&\leq&
\lim_{n\to \infty}\Pb\left(\sum_{i=1}^n(Y^{\prime}_{i,n}-p^{\prime}_{n,4r})\leq - a^{\prime}_n\right)\\ &\leq&
\lim_{n\to \infty}\frac{n p^{\prime}_{n,4r}}{(a^{\prime}_n)^2}=0.
\end{eqnarray*}
This concludes the proof of Theorem 4.1. \done

\medskip
\noindent {\textit{Proof of Theorem 4.2.}} With the  noncentrality parameter $\delta_n$ chosen as $\delta_n=2r\log n$, $0<r<1$,
one has the same tail behavior as in the normal case. Namely (see Section 5 of \cite{DJ} for details),
\begin{eqnarray*}
\Pb(\chi^2_\nu(0)>2s\log n)&=&O\left( n^{-s}\log^{-1/2}n\right),\quad  0<s\leq 1,\\
\Pb(\chi^2_\nu(\delta_n)>2s\log n) &=&O\left( n^{-(\sqrt{s}-\sqrt{r})^2}\log^{-1/2}n\right),\quad 0<r<s\leq 1.
\end{eqnarray*}
 With these relations available,
one immediately gets the analog of (\ref{asympt}) for the model in hand, and the analysis  proceeds exactly as
in the normal case. \done

\section{Concluding remarks}

\noindent
In this paper we study a new family of goodness-of-fit test statistics
that have the form of the empirical process in weighted sup-norm metrics with
EKFP weight functions $q$. We call these EDF-based statistics
the Cs\"{o}rg\H{o}-Cs\"{o}rg\H{o}-Horv\'{a}th-Mason (CsCsHM) statistics, after M. Cs\"{o}rg\H{o}, S. Cs\"{o}rg\H{o}, Horv\'{a}th, and Mason
whose beautiful results in \cite{CsCsHM} were a starting point of the current study.
 Taking $q$ as in (\ref{qq}), we
compare the corresponding one-sided statistic to the higher criticism statistic ${\rm HC}^+_n$ which is normalized by the SDP weight function $q(u)=\sqrt{u(1-u)}$.
Since, under $H_0$, the statistic ${\rm HC}^+_n$ tends  to $\infty$ in probability, and even almost surely,  the problem of determining the critical value
of the corresponding test procedure is unavoidable.
A resolution of this problem was guessed in Section 3 of \cite{DJ} via the Eicker-Jaeschke limit theorem for an
accordingly normalized ${\rm HC}^+_n$ sequence of statistics as in (\ref{DJ}). A correct version along these lines obtained in our Proposition 2.1,
concludes, regardless of a particular value of $0<\alpha_0<1$, an extreme value distribution that does
not depend on the parameter $\alpha_0$ at all. Consequently, it does not provide an appropriate limit distribution for any use of the ${\rm HC}^+_n$
statistics, theoretical and practical alike.

 In this paper, we have resolved this problem by using an entirely different strategy based on the theory of weighted empirical processes.
 In particular, we have shown that, as long as one deals with sup-norm functionals of weighted empirical processes, the most natural family of weights consists of
 the EFKP upper-class functions of a Brownian bridge.
 An immediate advantage of our approach is the appropriate identification of the limit distributions of the test statistics in hand
 under the null hypothesis. By using the algorithm in Section 6, these limit distributions are tabulated. Numerical comparison
 of the CsCsHM confidence bands obtained in Section 3.2 with the corresponding Kolmogorov-Smirnov confidence bands suggests that
 the CDF's of the CsCsHM statistics may
 converge to their respective limit CDF's at the rate of $O(n^{-1/2})$, the same rate as that in (\ref{rate}).
 In addition, we have shown that, like the higher criticism procedure, the whole class of test statistics $T_n^+(q,I)$, $I=(a,b)\subseteq(0,1)$, has the optimal adaptivity property
 (see Theorems 4.1--4.2 and Remark 4.1).

We also note that, when compared to the higher criticism statistic ${\rm HC}^+_n$, the CsCsHM test statistic (\ref{nasha})
provides a right solution in the sense that it does correctly the job that the former was intended to do
without requiring a large sample size like $n=10^6$ that, in case of the former, only indicated explosion
to infinity instead of slow convergence.
Therefore, in practical applications, rather than using the higher criticism statistic as in (\ref{DJT}), or as in (\ref{DJ2}), we recommend to use the test
statistic (\ref{nasha})
whose critical values are easily obtained by using the algorithm as in Section 6.

\section*{Acknowledgements} The authors wish to thank Miklos Cs\"{o}rg\H{o} for helpful discussions and suggestions.
The research of N. Stepanova was supported by an NSERC grant. The research of T. Pavlenko was partially supported by
the SRC grant C0595201 and by an NSERC grant during the author's stay at Carleton Univertsity in March 2014.

\noindent Natalia Stepanova, {School of Mathematics and Statistics, Carleton University, 1125 Colonel By Drive, Ottawa, ON K1S 5B6 Canada.}
E-mail: nstep@math.carleton.ca
\medskip

\noindent Tatjana Pavlenko, {Department of Mathematics, KTH Royal Institute of Technology, 106-91, Stockholm, Sweden.
E-mail: pavlenko@math.kth.se

\begin{thebibliography}{99}


\bibitem{ADarling} {Anderson, B. J.} \& {Darling, D. A.} (1952) Asymptotic theory of certain ``goodness of fit'' criteria based on stochastic processes.
\textit{Ann. Math. Statist.}, \textbf{23}, 193--212.

\bibitem{Bickel} {Bickel, P. J.} (1974) Edgeworth expansions in nonparametric statistics.
\textit{Ann. Statist.}, \textbf{2}, 1--20.

\bibitem{Borovkov} {Borovkov, A. A.} \& {Sycheva, N. M.} (1968) On some asymptotically optimal nonparametric tests. \textit{Teor. Verojatnost. i
Primenen.}, \textbf{13}, 385--418.

\bibitem{Cai} {Cai, T. T., Jeng, X. J,} \& {Jin, J.} (2011) Optimal detection of heterogeneous and heteroscedastic
mixtures. \textit{J. R. Statist. Soc. B}, \textbf{73}, 629--662.

\bibitem{Chib} {Chibisov, D. M.} (1964) Some theorems on the limiting behaviour of empirical
distribution functions. \textit{Selected Transl. Math. Statist. Probab.}, \textbf{6}, 147--156.

\bibitem{CsCsHM} {Cs\"{o}rg\H{o}, M., Cs\"{o}rg\H{o}, S., Horv\'{a}th, L.,} \& {Mason, D.} (1986) Weighted empirical and quantile processes.
\textit{Ann. Probab.}, \textbf{14}, 31--85.

\bibitem{CsHor} {Cs\"{o}rg\H{o}, M.} \& {Horv\'{a}th, L.} (1988) Nonparametric methods for changepoint problems. In:
\textit{Handbook of Statistics} \textbf{7} (eds. P. R. Krishnaiah \& C. R. Rao), 403--425. Elsivier Science Publishers B.V., Amsterdam.


\bibitem{CsH} {Cs\"{o}rg\H{o}, M.,} \& {Horv\'{a}th, L.} (1993) \textit{Weighted approximations in probability and statistics}.
 Wiley, New York.

\bibitem{CsHorvath} {Cs\"{o}rg\H{o}, M.} \& {Horv\'{a}th, L.} (1997) \textit{Limit theorems in change-point analysis.}
 Wiley, New York.

\bibitem{CsR} {Cs\"{o}rg\H{o}, M.} \& {R\'{e}v\'{e}sz, P.} (1981) \textit{Strong approximation in probability and statistics}.
Academic Press, New York (Akad\'{e}miai Kiad\'{o}, Busdapest).


\bibitem{Darling} {Darling, D. A.} \& {Erd\H{o}s, P.} (1956) A limit theorem for the maximum of normalized sums of
independent random variables. \textit{Duke Math. J.}, \textbf{23}, 143--145.

\bibitem{DasG} {DasGupta, A.} (2008) \textit{Asymptotic Theory of Statistics and Probability}. Springer Science+Business Media, LLC.

\bibitem{DJ} {Donoho, D.} and {Jin, J.} (2004) Higher criticism for detecting sparse heterogeneous mixtures. \textit{Ann. Statist.},
\textbf{32}, 962--994.


\bibitem{DJ2009} {Donoho, D.} \& \textsc{Jin, J.} (2009) Feature selection by higher criticism thresholding achieves the optimal phase diagram.
\textit{Phil. Trans. R. Soc. A}, \textbf{367}, 4449--4470.

\bibitem{EE} {Eastwood, B. J.} \& {Eastwood, V. R.} (1998) Tabulating weighted functionals of Brownian bridges
via Monte Carlo simulation. In: \textit{Asymptotic methods in probability and statistics --- A volume in honour of Mikl\'{o}s Cs\"{o}rg\H{o}}
(ed. B. Szyszkowicz), 707--719. Elsivier Science B.V., Amsterdam.

\bibitem{Eicker} {Eicker, F.} (1979) The asymptotic distribution of the suprema of the standardized empirical processes. \textit{Ann. Statist.},
\textbf{7}, 116--138.

\bibitem{Gnedenko} {Gnedenko, B. V., Korolyuk, V. S.,} \& {Skorohod, A. V.} (1960) Asymptotic expansions in probability theory.
In: \textit{Proc. Fourth Berkley Symp. Math. Statist. Prob.}, \textbf{2}, 153--169, Univ. California Press.

\bibitem{Ing97} {Ingster, Yu. I.} (1997) Some problems of hypothesis testing leading to infinitely divisible distribution.
\textit{Math. Meth. Statist.}, \textbf{6}, 47--69.

\bibitem{Ing99} {Ingster, Yu. I.} (1998) Minimax detection of a signal for $l_n$-balls.
\textit{Math. Meth. Statist.}, \textbf{7}, 401--428.

\bibitem{IPT} {Ingster, Yu. I., Pouet, C.,} \& {Tsybakov, A. B.} (2009) Classification of sparse high-dimensional vectors.
\textit{Phil. Trans. R. Soc. A}, \textbf{367}, 4427--4448.

\bibitem{J} {Jaeschke, D.} (1979) The asymptotic distribution of the supremum of the standardized empirical distribution function on subintervals.
\textit{Ann. Statist.}, \textbf{7}, 108--115.

\bibitem{Jager} {Jager, L.} \& {Wellner, J. A.} (2007) Goodness-of-fit test via phi-divergences. \textit{Ann. Statist.},
\textbf{35}, 2018--2035.

\bibitem{Klaus} {Klaus, B.} \& {Strimmer, K.} (2013) Signal detection for rare and weak features: higher criticism or
false discovery rates? \textit{Biometrics}, \textbf{14}, 129--143.

\bibitem{MR} {Meinshausen, N.} \& {Rice, J.} (2006) Estimating the proportion of false null hypotheses among
a large number of independently tested hypotheses. \textit{Ann. Statist.}, \textbf{34}, 373--393.

\bibitem{Nikitin} {Nikitin, Ya. Yu.} (1995) \textit{Asymptotic efficiency of nonparametric tests}. Cambridge Univ. Press.

\bibitem{OP} {Orasch, M.} and {Pouliot, W.} (2004) Tabulating weighted sup-norm functionals used in change-point problem. \textit{J. Stat. Comput. Simul.},
\textbf{74}, 249--276.

\bibitem{OR} {O'Reilly, N.} (1974) On the weak convergence of empirical processes in sup-norm metrics.
\textit{Ann. Probab.}, \textbf{2}, 642--651.

\bibitem{Pavlenko12} {Pavlenko, T., Bj\"orkstr\"om, A.,} \& {Tillander, A.} (2012)  Covariance structure approximation via gLasso
in high-dimensional supervised classification. \textit{J. Appl. Stat.}, \textbf{8}, 1643--1666.

\bibitem{Shorak} {Shorack, G. R.} \& {Wellner, J. A.} (1986) \textit{Empirical processes with applications to statistics}. Wiley, New York.

\bibitem{Tillander13} {Tillander, A.} (2013) {Classification models for high-dimensional data with sparsity patterns}. PhD dissertation.  Stockholm University, Stockholm.


\end{thebibliography}
\end{document}